\documentclass[11pt]{amsart}
\usepackage{amssymb}
\usepackage{amsmath}
\usepackage{amsthm}
\usepackage[textwidth=16cm, textheight=24.5cm, marginratio=1:1]{geometry}

\usepackage{newtxtext} 
\usepackage[english]{babel}
\usepackage[utf8]{inputenc}
\usepackage[T1]{fontenc}
\usepackage{todonotes}
\usepackage{xcolor}
\usepackage[backref=page]{hyperref}
\hypersetup{
    colorlinks,
    linkcolor={red!80!black},
    citecolor={blue!80!black},
    urlcolor={blue!80!black}
}
\usepackage{setspace}
\usetikzlibrary{matrix,arrows, cd, calc}  
\newsavebox{\pullback}
\sbox\pullback{%
\begin{tikzpicture}%
\draw (0,0) -- (1ex,0ex);%
\draw (1ex,0ex) -- (1ex,1ex);%
\end{tikzpicture}}
\usepackage{enumitem}
\usepackage{booktabs}
\numberwithin{equation}{section}
\newtheorem{theorem}[equation]{Theorem}

\newtheorem{lemma}[equation]{Lemma}

\newtheorem{proposition}[equation]{Proposition}
\theoremstyle{definition}
\newtheorem{definition}[equation]{Definition}

\newtheorem{remark}[equation]{Remark}
\newtheorem{example}[equation]{Example}

\DeclareMathOperator{\coker}{coker}
\DeclareMathOperator{\im}{im}
\DeclareMathOperator{\Spec}{Spec}

\DeclareMathOperator{\Alg}{Alg}
\DeclareMathOperator{\CQAlg}{CQAlg}
\DeclareMathOperator{\Coh}{Coh}
\DeclareMathOperator{\CQCoh}{CQCoh}
\DeclareMathOperator{\pr}{pr}

\newcommand{\onto}{\twoheadrightarrow}
\newcommand{\into}{\hookrightarrow}

\newcommand{\kk}{\Bbbk}%
\newcommand{\spann}[1]{\left\langle #1 \right\rangle}
\DeclareMathOperator{\Ann}{Ann}
\DeclareMathOperator{\OHilb}{Hilb}
\DeclareMathOperator{\Quot}{Quot}
\DeclareMathOperator{\CQuot}{CQuot}
\newcommand{\OO}{\mathcal{O}}%
\newcommand{\OOhat}{\widehat{\OO}}%
\DeclareMathOperator{\GL}{GL}%
\DeclareMathOperator{\Isom}{Isom}
\DeclareMathOperator{\tr}{tr}%
\DeclareMathOperator{\id}{id}%
\DeclareMathOperator{\gl}{\mathfrak{gl}}%
\DeclareMathOperator{\sll}{\mathfrak{sl}}%
\DeclareMathOperator{\PGL}{PGL}%
\DeclareMathOperator{\End}{End}%
\DeclareMathOperator{\Sym}{Sym}%
\newcommand{\Gmult}{\mathbb{G}_{\mathrm{m}}}
\DeclareMathOperator{\Iar}{Iar}%
\DeclareMathOperator{\Gr}{Gr}%
\newcommand{\mm}{\mathfrak{m}}%
\newcommand{\gr}{\operatorname{gr}}%
\newcommand{\orr}{\mathrm{or}}
\DeclareMathOperator{\Gor}{Gor}
\newcommand{\sym}{\mathrm{sym}}
\DeclareMathOperator{\CQ}{CQ}
\newcommand{\cE}{\mathcal{E}}
\DeclareMathOperator{\CCQ}{CCQ}
\newcommand{\BBname}{Bia{\l}ynicki-Birula}%

\DeclareMathOperator{\rk}{rk}
\DeclareMathOperator{\sm}{sm}

\newcommand{\kkt}{\kk[\![t]\!]}
\newcommand{\kkbar}{\overline{\kk}}%

\newcommand{\Xbar}{\overline{X}}%

\begin{document}

\title{The Iarrobino scheme: a self-dual analogue of the Hilbert scheme of
points}
\author{Joachim Jelisiejew}
\address[J. Jelisiejew]{University of Warsaw, Faculty of Mathematics,
Informatics and Mechanics, Banacha 2, 02-097 Warsaw, Poland}
\email{j.jelisiejew@uw.edu.pl}
\thanks{Supported by National Science Centre grant 2023/50/E/ST1/0033.}
\date{\today{}}

\begin{abstract}
For a fixed quasi-projective scheme $X$ we introduce a self-dual analogue of
$\OHilb_d(X)$ which we call the Iarrobino scheme of $X$. It is a fine
moduli space of zero-dimensional subschemes of $X$ together with a self-dual filtration. It contains an open subset
that parameterises oriented Gorenstein subschemes.
    We provide also self-dual analogues of the Quot scheme of points and of the
    stacks of coherent sheaves and finite algebras. We prove that the resulting analogues of Hilbert and
    Quot schemes are smooth for $X$ a smooth curve.
    The Iarrobino scheme is a natural ambient space for enumerative invariants of algebras as introduced by Micha{\l}ek.
    When $X$ is an affine space, it also admits a commuting matrices description and related virtual structure.
\end{abstract}
\maketitle

\section{Introduction}

\newcommand{\fullrkarg}[1]{\mathbb{P}(\Sym^2#1)^{\mathrm{full\,rk}}}
\newcommand{\fullrk}{\fullrkarg{V}}
\newcommand{\HilbdX}{\OHilb_d(X)}%
\newcommand{\HilbdGororX}{\OHilb_d^{\Gor,\orr}(X)}%
\newcommand{\HilbdGororsmX}{\OHilb_d^{\Gor,\orr, \circ}(X)}%
\newcommand{\HilbdGorX}{\OHilb_d^{\Gor}(X)}%
\newcommand{\HilbdGororAn}{\OHilb_d^{\Gor,\orr}(\mathbb{A}^n)}%
\newcommand{\HilbdAn}{\OHilb_d(\mathbb{A}^n)}%
\newcommand{\Oext}{\widetilde{O}}%

    Let $\HilbdX$ be the Hilbert scheme of points on a quasi-projective $\kk$-scheme
    $X$. The $\kk$-points of $\HilbdX$ are finite degree $d$ closed subschemes
    $Z\subseteq X$. For $X$ a smooth curve or a smooth surface, the
    Hilbert scheme $\HilbdX$ is smooth.

    In this work we introduce the Iarrobino scheme of points $\Iar_d(X)$,
    which serves as analogue of $\HilbdX$ in the self-dual setting.
    The scheme $\Iar_d(X)$ parameterises pairs $(Z\subseteq X,\, [q_{\bullet}])$, where
    $Z$ is as above and $q_{\bullet}$ is a \emph{broken} quadric (see~\eqref{eq:kpointsCQ}) that makes
    \emph{an
    associated graded of} the sheaf $\OO_Z$ self-dual. An important case is
    when $q_{\bullet}$ is a single full rank quadric and $\OO_Z$ itself is self-dual; this
    is the oriented Gorenstein locus. The forgetful map
    \[
        \Iar_d(X)\to \HilbdX
    \]
    is projective.
    We discuss $\Iar_d(X)$ more precisely in~\S\ref{secintro:IarrobinoScheme} below.

    The motivation for $\Iar_d(X)$ comes from the following
    sources:
    \begin{enumerate}[topsep=0pt,itemsep=0pt]
        \item As an ambient scheme for the Iarrobino's symmetric decomposition
            of the Hilbert function~\cite{iarrobino_associated_graded}, which
            is a central tool for classifying and
            understanding finite non-graded Gorenstein
            algebras.
        \item As an ambient space for enumerative invariants of tensors~\S\ref{ssec:junehuh}, as developed by
            Mateusz Micha{\l}ek and others~\cite{Michalek__survey,
            Jagiella_Pielasa_Shatshila} following the work of June Huh on matroids.
        \item As a projective scheme containing of the locus of oriented
            Gorenstein subschemes and preserving its self-duality.
        \item As a symmetric version of ADHM construction of $\OHilb_d(\mathbb{A}^n)$ via commuting matrices.
            In particular this equips $\Iar_d(\mathbb{A}^2)$ with a virtual fundamental class, see~\S\ref{ssec:introADHM} below.
    \end{enumerate}
    The construction of the Iarrobino scheme has two principal ingredients: Hilbert schemes and 
    the variety of completed quadrics, which dates back to classical
            times~\cite{Kleiman__history, Semple, Tyrrell__complete_quadrics},
            and is reinterpreted as a wonderful
            compactification~\cite{Concini_Procesi__Complete_symmetric_varieties}.

    \subsection{Completed quadrics and the Iarrobino
    scheme}\label{secintro:IarrobinoScheme}

    Throughout, $\kk$ is a field of characteristic different from $2$.
    A finite $\kk$-scheme $Z$ is \emph{Gorenstein} if its dualising sheaf
    $\omega_Z = \OO_Z^{\vee}$, see~\S\ref{ssec:duality}, is invertible, or,
    equivalently, trivial: there is an isomorphism $q\colon \omega_Z\to
    \OO_Z$. For a variety $X$, the locus of Gorenstein schemes $\HilbdGorX\into \HilbdX$ is
    open~\cite[p.2056]{casnati_notari_irreducibility_Gorenstein_degree_9},
    but it is usually far from being projective even if $X$ is projective, see
    Example~\ref{ex:associatedGradedSymmetric}. The
    Iarrobino scheme remedies this.

    \begin{theorem}[{Definition~\ref{ref:Iar:def}, Proposition~\ref{ref:Iarkpoints:prop}}]\label{refintro:existence:thm}
        Let $X$ be a quasi-projective $\kk$-scheme.
        There is a quasi-projective scheme $\Iar_d(X)$ with a
        \emph{projective} morphism
        \[
            \tau_X\colon \Iar_d(X) \to \HilbdX,
        \]
        such that the $\kk$-points of $\Iar_d(X)$
        correspond bijectively to sequences (of varying length):
        \[
            X\supseteq Z_0\supsetneq Z_1\supsetneq Z_2\supsetneq \ldots
        \]
        of zero-dimensional subschemes of $X$ such that $\deg Z_0 = d$,
        together with a choice (up to scalar multiplication) of $\OO_X$-module isomorphisms
        \begin{equation}\label{eq:isos}
            q_i\colon \left(\frac{I_{Z_{i+1}}}{I_{Z_i}}\right)^{\vee}\to
            \frac{I_{Z_{i+1}}}{I_{Z_i}},
        \end{equation}
        for $i=0,1, \ldots $, which satisfy $q_i^{\top} = q_i$. The map
        $\tau_X$ sends such a $\kk$-point to $[Z_0]\in \HilbdX$.
    \end{theorem}
    Less formally put, the scheme $Z_0$ comes with an exhaustive filtration by
    self-dual quotients.

    The ``good'' locus in Theorem~\ref{refintro:existence:thm} corresponds to
    the case with
    only $Z_0$. Here, the subquotient $I_{Z_{1}}/I_{Z_0}$ is isomorphic to $\OO_{Z_0}$ and $q_0\colon
    \OO_{Z_0}^{\vee} \to \OO_{Z_0}$ is the trivialisation of the dualising
    sheaf, hence $Z_0$ is Gorenstein. This locus is open in
    $\Iar_d(X)$, see Definition~\ref{ref:unbrokenopen:lem}. The cases with
    ``residual'' quadrics $q_1$, $q_2$, etc. serve to compactify it.

    The Iarrobino scheme, just like the usual Hilbert scheme of points, comes
    accompanied by the Quot scheme version: for a coherent sheaf $\cE$ on $X$,
    we have the \emph{self-dual Quot
    scheme of points} $\CQuot_d(\cE)$, see Definition~\ref{ref:Quot:def}. As
    with Hilb and Quot, we have $\Iar_d(X) \simeq \CQuot_d(\OO_X^{\oplus 1})$.
    We also obtain the stacky counterparts: the self-dual Artin stack $\CQAlg_d$ which plays the
    role of  the stack $\Alg_d$ of finite algebras, and the self-dual Coherent
    Stack $\CQCoh_d X$ on a quasi-projective scheme $X$, which plays the role
    of the stack of coherent sheaves $\Coh_d X$, for both
    see~\S\ref{ssec:stacks}.
    The main emphasis of this paper is on constructing the above schemes.
    Additionally, we provide the following results in the case of curves.
    \begin{theorem}[Theorem~\ref{ref:smoothnessCQuotCurve:thm}]\label{refintro:smoothForCurves:thm}
        Let $C$ be a smooth irreducible curve and $\cE$ be a locally free sheaf on $C$.
        Then the scheme $\CQuot_d(\cE)$ is smooth connected of dimension
        $d\cdot (1+\rk \cE) - 1$.
    \end{theorem}
    \begin{theorem}[Theorem~\ref{ref:smoothnessCurve:thm}]\label{refintro:IarrobinoToHilbert:thm}
        Let $C$ be a smooth connected curve.
        The scheme $\Iar_d(C)$ is smooth irreducible of dimension
        $2d-1$. The Iarrobino-to-Hilbert morphism
        \[
            \tau_C\colon \Iar_d(C)\to \OHilb_d(C)
        \] is flat and projective with integral fibres of dimension $d-1$
        which are reduced complete intersections in $\Iar_d(C)$. The morphism satisfies
        $(\tau_C)_* \OO_{\Iar_d(C)} = \OO_{\OHilb_d(C)}$.
    \end{theorem}

    \newcommand{\QuotX}[2]{\Quot_{#1}\left( \OO_X^{\oplus #2} \right)}%
    \newcommand{\CQuotX}[2]{\CQuot_{#1}\left( \OO_X^{\oplus #2} \right)}%
    \newcommand{\CQuotpar}[3]{\CQuot_{#1}\left( \OO_{#3}^{\oplus #2} \right)}%

    The Hilbert scheme of points admits a fairly elementary
    construction~\cite{Haiman__exposition, Huibregtse_construction,
    Gustavsen_Laksov_Skjelnes__Elementary_explicit_const},
    \cite[Chapter~18]{Miller_Sturmfels}.  The construction of the Iarrobino
    scheme involves one
    additional element: the variety of
    completed quadrics.
    This variety itself appears also as a special case of the construction of the Quot version.
    In the classical setup, the Quot scheme
    $\Quot_{d}\left(\OO_{\Spec(\kk)}^{\oplus d}\right) = \Gr(d, d)$ is a
    single point
    corresponding to the $d$-dimensional vector space $V = \kk^{\oplus
    d}$. In contrast, the
    self-dual Quot scheme $\CQuotpar{d}{d}{\Spec(\kk)}$ is
    isomorphic to the \emph{variety of completed quadrics}
    $\CQ(V)$.
    This is a smooth projective variety of dimension $\binom{d+1}{2}-1$, which can be defined as the closure
    of the embedding
    \begin{equation}\label{eq:introproduct}
        \fullrk \into \mathbb{P}(\Sym^2 V) \times \mathbb{P}\left(
                \Sym^2\Lambda^2 V \right)\times
                \ldots \times \mathbb{P}\left( \Sym^{2}\Lambda^{d-1}V
                \right) \times
                \mathbb{P}\left(\Sym^2\Lambda^{d} V\right)
    \end{equation}
    given by the exterior powers. We discuss it
    in~\S\ref{sec:completedQuadrics}.

    \goodbreak
    \subsection{Applications}\label{secintro:motivation}

        \subsubsection{Characteristic numbers of algebras}\label{ssec:junehuh}

            June Huh and others analysed realisable
            matroids by considering the closure of its span in the
            permutohedral toric variety.
            The permutohedral toric variety can be realised as the closure of
            a torus inside the variety of completed quadrics $\CQ(V)$, see
            also Example~\ref{ex:completedQuadricsOnSmooth}.
            Conner and Micha{\l}ek~\cite{Conner_Michalek} generalised this and proposed to study characteristic numbers of
            algebras (and, more generally, tensors) via closures in $\CQ(V)$ itself.

            In this setup, $V$ is an algebra and the multiplication map
            $\Sym^2
            V\to V$ yields a subspace $i\colon V^{\vee} \into \Sym^2 V^{\vee}$. If $V$ is
            Gorenstein, the image of $i$ contains full rank quadrics. One
            takes the closure of $(\im i) \cap \fullrk$ inside $\CQ(V^{\vee})
            \simeq \CQ(V)$ and obtains a subvariety
            \[
                X_V \subseteq \CQ(V)
            \]
            of
            dimension $\dim_{\kk} V - 1 = d-1$.
            The variety $\CQ(V)$ has
            Picard group
            \[
                \mathbb{Z}L_1 \oplus \ldots \oplus \mathbb{Z}L_{d-1}
            \]
            for certain distinguished line bundles $L_1, \ldots
            ,L_{d-1}$, which are pullbacks of $\OO(1)$-bundles from respective
            factors of~\eqref{eq:introproduct}.
            The above mentioned \emph{characteristic number} of the algebra $V$ is any intersection number
            $L_1^{b_1} \ldots L_{d-1}^{b_{d-1}}|_{X_V}$ where $b_1, \ldots
            ,b_{d-1}\in \mathbb{Z}_{\geq 0}$ sum up to $d - 1$.

            \medskip
            Our contribution here is as follows.
            In general, the variety $X_V$ is hard to understand due to the closure and the varieties related to various $V$ are a priori unrelated.
            We provide (see~\S\ref{ssec:charNumbers}) an interpretation of $X_V$ as the
            closure of the Gorenstein part of the fiber over $[\Spec(V)]$ of the
            Iarrobino-to-Hilbert map.

            For a smooth curve, the Iarrobino-to-Hilbert map $\tau$ is flat
            with irreducible fibers, see
            Theorem~\ref{refintro:IarrobinoToHilbert:thm} above, so the
            varieties $X_V$ coincide with the fibers
            and we obtain
            immediately the following result, first proven by
            Jakub Jagie{\l}{\l}a, Pawe{\l} Pielasa, and Anatoli Shatsila by a direct calculation.
            \begin{theorem}[{\cite[Theorem~1.3, Corollary~1.4]{Jagiella_Pielasa_Shatshila}},
            see Proposition~\ref{ref:characteristicNumbers:prop}]\label{ref:Jagiella:introthm}
                For every finite degree $d$ algebra of the form $\kk[x]/(f)$,
                the characteristic numbers are the same.
            \end{theorem}

        \subsubsection{Classification of zero-dimension Gorenstein
            algebras and origin of the name}\label{ssec:IarrobinoIntro}

    \newcommand{\annmm}[1]{(0:\mm^{#1})}
    \newcommand{\DeltaHilb}[1]{\Delta_{#1}}%
    \newcommand{\QHilb}[1]{\mathcal{Q}(#1)}%
            The classification of finite degree $d$ algebras is a classical
            subject. The theory reduces immediately to the case of local
            $\kk$-algebras $(A, \mm, \kk)$.
            A key numerical invariant of such an algebra $A$ is its Hilbert
            function $H_A(i) := \dim_{\kk} \mm^i/\mm^{i+1}$ obtained from its
            associated graded algebra
            \[
                \gr(A) = \bigoplus_{i=0}^{\infty} \frac{\mm^i}{\mm^{i+1}}.
            \]
            Since $A$ is Artinian, the sum above is actually finite, ending at
            the
            maximal number $s$ such that $\mm^s \neq 0$.
            For Gorenstein algebras, their classification and
            numerical invariants are also very actively researched, see, for
            example~\cite{elias_rossi_short_Gorenstein,
            Jelisiejew_Masuti_Rossi, Iarrobino_Macias_Marques}.
            The most important numerical invariant, by Anthony
            Iarrobino~\cite{iarrobino_associated_graded}, is the
            \emph{symmetric decomposition} of $H_A$:
            \[
                H_A = \sum_{\delta=0}^{s-2} \Delta_{\delta}
            \]
            where the $s-1$ functions $\Delta_{0}, \ldots , \Delta_{s-2}$ such that
            $\Delta_{\delta}(i) = \Delta_{\delta}(s-i-\delta)$ for every
            $0\leq i\leq s - \delta$ and $\Delta_{\delta}(i) = 0$ otherwise.
            The symmetric decomposition refines $H_A$ and restricts possible
            Hilbert functions of Gorenstein algebras, for example
            $H_A =(1,n,n+1,1,0,0, \ldots )$ is not possible for any $n$.
            Still, possible $H_A$ and $\Delta_{\bullet}$ are far from being
            classified.

            The symmetric decomposition is induced from a filtration on
            $\gr(A)$ obtained from the Loewy filtration $(0:\mm^{s+1-i})_{i\geq 0}$
            on $A$. The associated graded of this filtration is a
            self-dual $\gr(A)$-module
            $\bigoplus_{\delta=0}^{s-2} \QHilb{\delta}$,
            see~\S\ref{ssec:symmDec} for details.

            Let us pass from numerics to geometry.
            The associated graded $\gr(A)$ can be
            realised by embedding $\Spec(A) \subseteq \mathbb{A}^{n}$ with
            support at the origin and as taking the torus limit at $t\to
            \infty$ for the usual action $t\cdot (x_1, \ldots ,x_n) =
            (tx_1, \ldots ,tx_n)$, see~\S\ref{ssec:IarrobinoCurves}. This simple
            observation leads, via the \BBname{} decomposition, to significant
            results in the geometry of $\OHilb_d(\mathbb{A}^n)$, see e.g.~\cite{Jelisiejew__Elementary}.

            In this work we provide a Gorenstein analogue of the above. We realise $\bigoplus_{\delta=0}^{s-2}
            \QHilb{\delta}$ as a limit of a Gorenstein $\Spec(A)\subseteq
            \mathbb{A}^{n}$ in the \emph{Iarrobino
            scheme}, see Theorem~\ref{ref:name:thm} and the picture below.

            \begin{center}
            \resizebox{0.75\textwidth}{!}{
                \begin{tikzpicture}
                    \node (iardesc) at (-10.3, 2.5) {\large $\lim_{t\to \infty} t\cdot
                    (A,[q_0]) = (\QHilb{\bullet}, [q_{\bullet}])$};
                    \node (hilbdesc) at (-10.3, -0.5) {\large $\lim_{t\to \infty} t\cdot
                        [A] = [\gr(A)]$};
                    \node (iar) at (-5.3, 2.5) {{\Large $\operatorname{Iar}_d(X)$}};
                    \node (hilb) at (-5.3, -0.5) {{\Large $\operatorname{Hilb}_d(X)$}};
                    \path[->,line width=1pt] (iar) edge (hilb);
                    \node (tau) at (-5.8, 1.0) {\Large $\tau_{X}$};
                    \node (orbits) [blue] at (-2, 1.0) {\large
                        $\mathbb{G}_{\mathrm{m}}$-orbit closures};

                        \begin{scope}[shift={(0, -3)}]
                            \draw[blue] plot [smooth] coordinates {(-3, 2.3) (-2.5, 2.2) (-2, 2.5)
                            (-1.1, 2.5) (-1, 2.5)};
                        \end{scope}
                        \begin{scope}[shift={(0, 3.2)}]
                            \draw[blue] plot [smooth] coordinates {(-3, -1) (-2.5, -1.1) (-2, -0.8)
                            (-1.1, -0.7) (-1, -0.3)};
                        \end{scope}
                        \draw plot [smooth cycle] coordinates {(0, 2.5) (-1, 1.5) (-2, 1.85)
                        (-3.5, 2) (-3, 3) (-2, 2.75) (-0.7, 3.2)};
                        \draw plot [smooth cycle] coordinates {(0, -0.3) (-0.2, -0.9) (-0.9, -1) (-1.3, -1.4)
                            (-1.5, -1.6) (-2.5, -1.6) (-3.2, -1.7) (-4, -1) (-3.65, -0.5) (-3.8, -0.1)
                        (-3, 0) (-1.5, 0.3) (-0.7, 0)};
                    \end{tikzpicture}}
                \end{center}
                This allows to geometrically understand the symmetric decomposition via \BBname{} decompositions on
            $\Iar_d(\mathbb{A}^n)$.

            \subsubsection{Commuting matrices}\label{ssec:introADHM}
        \newcommand{\comp}{\mathcal{C}}%
        \newcommand{\anticomp}{a\mathcal{C}}%

                In the case $X = \mathbb{A}^n$, the Hilbert scheme of points $d$
                points admits a description known as the ADHM
                construction~\cite[\S2]{nakajima_lectures_on_Hilbert_schemes}.
                Fix a $d$-dimensional vector space $V$. Then $\HilbdAn$ is the free quotient of
                \[
                    \left\{ (X_1, \ldots ,X_n, v)\in \End(V)^{\oplus n}\times V\ |\
                    \kk[X_1, \ldots ,X_n]v = V\ \mbox{and}\ \forall_{i,j}\ X_iX_j = X_jX_i \right\}
                \]
                by the natural $\GL(V)$-action. This description is quite useful in
                enumerative geometry, yielding the global Kuranishi chart for $n=2$,
                and in moduli spaces, as it tightly links the Hilbert and Quot schemes
                to the variety of commuting
                matrices~\cite{jelisiejew_sivic}.

                It is natural to ask what happens when we replace $\End(V)$ by
                some classical subalgebra. Here we discuss the orthogonal group.
                More precisely, let
                $V$ come with a fixed full rank quadric $Q$ and let $\Oext(V,
                Q)\subseteq \GL(V)$ be
                the stabiliser of $[Q]\in \mathbb{P}\Sym^2(V)$, so that $\Oext(V, Q) =
                O(V, Q)\cdot \Gmult$, where $O(V, Q)$ is the usual orthogonal group
                and $\Gmult\subseteq \GL(V)$ are the scalar matrices. The quadric $Q$
                identifies $V$ with $V^{\vee}$ and in particular allows us to speak
                about symmetric matrices $\End^{\sym}(V)\subseteq \End(V)$. We
                have the following analogue:

                \begin{proposition}\label{ref:ADHMGorenstein:prop}
                    The oriented Gorenstein locus $\HilbdGororAn$ is a free quotient
                    of
                    \[
                        \left\{ (X_1, \ldots ,X_n, v)\in \End^{\sym}(V)^{\oplus n}\times V\ |\
                        \kk[X_1, \ldots ,X_n]v = V\ \mbox{and}\ \forall_{i,j}\ X_iX_j = X_jX_i \right\}
                    \]
                    by the natural $\Oext(V, Q)$-action. Equivalently, it is a free quotient by $\GL(V)$-action of
                    the locus
                    \[
                        \left\{ ([Q], X_1, \ldots ,X_n, v)\in \fullrk \times \End^{\sym, Q}(V)^{\oplus n}\times V\ |\
                        \kk[X_1, \ldots ,X_n]v = V\ \mbox{and}\ \forall_{i,j}\ X_iX_j = X_jX_i \right\}
                    \]
                    where $\End^{\sym, Q}(V)$ means endomorphisms symmetric with respect to $Q$.
                \end{proposition}
                In the above language, the lack of properness of $\HilbdGororAn\to\Sym^d \mathbb{A}^n$ can be traced
                to lack of properness of $\fullrk$.
                We address it by compactifying $\fullrk$ to $\CQ(V)$.
                One delicate point is how to extend $\End^{\sym, Q}(V)$ outside $\fullrk$. We show (\S\ref{ssec:compatibleBundle}) that indeed there exists
                a subbundle $\comp\into \End(V)\times \CQ(V)$ such that its fibre over every $[Q]\in \fullrk \subseteq \CQ(V)$ is equal to $\End^{\sym, Q}(V)$.
                There exists also the \emph{anti}symmetric analogue which we call the anticompatible bundle $\anticomp \to \CQ(V)$.

                \begin{theorem}\label{ref:ADHMdesciption:thm}
                    The scheme $\Iar_d(\mathbb{A}^n)$ is a free quotient by $\GL(V)$ of the space
                    \[
                        \left\{ (X_1, \ldots ,X_n, v)\in \comp^{\oplus n}\times V\ |\
                        \kk[X_1, \ldots ,X_n]v = V\ \mbox{and}\ \forall_{i,j}\ X_iX_j = X_jX_i \right\}
                    \]
                    where $\comp^{\oplus n}\to \CQ(V)$ is the $n$-th direct sum of the compatible bundle.
                \end{theorem}
                Theorem~\ref{ref:ADHMdesciption:thm} is proven in~\S\ref{ssec:commuting}.
                One can introduce a smooth ambient space for $\Iar_d(\mathbb{A}^n)$ by taking the free $\GL(V)$-quotient of
                \[
                    \left\{ (X_1, \ldots ,X_n, v)\in \comp^{\oplus n}\times V\ |\
                    \kk[X_1, \ldots ,X_n]v = V\right\}.
                \]
                On this space, every commutator $[X_i, X_j]$ lies in the pullback of the anticompatible bundle $\anticomp$ that has rank $\binom{d}{2}$.
                In particular, $\Iar_d(\mathbb{A}^2)$ is the zero section of the pullback of $\anticomp$. This observation equips $\Iar_d(\mathbb{A}^2)$ with a
                virtual fundamental class in dimension
                \[
                    \binom{d+1}{2}-1+d+2\binom{d+1}{2} - d^2 - \binom{d}{2} = 3d-1,
                \]
                which is exactly the dimension of the oriented Gorenstein locus.

        \subsubsection{Enumerative geometry of the Iarrobino scheme}

            Enumerative geometry of Hilbert schemes in low dimensions is much investigated. We hope that the enumerative geometry of the Iarrobino scheme can reflect this.

            For a smooth curve $C$,
            the Iarrobino-to-Hilbert map $\tau =
            \tau_{\Iar_d(C)} \colon \Iar_d(C)\to \OHilb_d(C)$ is
            $\OO$-connected (see Theorem~\ref{ref:smoothnessCurve:thm}). But more is true and will be proven
            in subsequent work with Reinier F.~Schmiermann and Andrea Ricolfi. We decided to announce it here for 
            illustration. Namely,
            actually $R\tau_{*}\OO_{\Iar_{d}(C)} = \OO_{\OHilb_d(C)}$. Moreover, at least for $C = \mathbb{A}^1$, this equality
            generalises to Euler characteristics of the natural line bundles:
            \[
                \sum_{d\geq 0} \sum_{k=0}^d \chi(L_{k, d})q^dr^k = \left(
                \sum_{d\geq 0} \frac{1}{\prod_{i=1}^d(1-t^i)}q^d
                \right)\cdot \left( \sum_{d\geq 0}
                \frac{1}{\prod_{i=1}^d(1-t^i)} q^dr^d t^{(d-1)d} \right)
            \]
            where $t$ is the equivariant parameter for the usual $\Gmult$-action on $\mathbb{A}^1$,
            while for $i=1, \ldots ,d-1$
            the bundle $L_{i,d}$ on $\Iar_d(\mathbb{A}^1)$ is the pullback of $L_i$ discussed
            in relation to characteristic numbers in~\S\ref{ssec:junehuh}
            above, and $L_{0,d} = \OO_{\Iar_d(\mathbb{A}^1)}, L_{d, d}$ are defined likewise.

            For surfaces, one can hope for extending these results by using the perfect obstruction
            theory coming from the zero section presentation, as explained at the end of~\S\ref{ssec:introADHM} above.

\section*{Acknowledgements}

    The author is very grateful to Michel Brion, Maciej Do{\l}{\k{e}}ga, Martijn
    Kool, Mateusz Micha{\l}ek, Andrea Ricolfi, Reinier Schmiermann, and Jaros{\l}aw Wi{\'s}niewski for enlightening
    discussions. I additionally thank Andrea and Reinier for their insights
    and questions shared during our discussions. I thank Michel Brion also for
    directing me to~\cite{DeConcini_Springer}. Thanks to Alessio Sammartano
    and Richard Thomas for
    helpful comments.
    I thank Tony Iarrobino for permission to baptise $\Iar_d(X)$ using his
    name. Finally, thanks to the organisers of \emph{Hilbert scheme of points}
    day in Berlin, \emph{\textbf{G}eometry and \textbf{C}ommutative algebra}
    workshop in Milano, \emph{Enumerative geometry of the Hilbert scheme of
    points} in Les Diablerets, and the MFO
    workshop \emph{Recent Trends in Algebraic Geometry} for giving me the
    possibility to present this work.

\section{Preliminaries}

We work over a base field $\kk$. We assume that $2\in \kk$ is
invertible. For a $\kk$-vector space $V$, we denote by
$V^{\vee}$ the dual space.

\subsection{Duality for zero-dimensional modules}\label{ssec:duality}
In this section we collect basic facts about duality.
A $\kk$-algebra is \emph{finite} if it is finite-dimensional over $\kk$.
  If $S$ is a $\kk$-algebra and $M$ is an
finite $S$-module, then the dual vector space $M^{\vee}$ is also an $S$-module via contraction:
$(s\cdot \varphi)(m) := \varphi(sm)$
for all $s\in S$, $\varphi\in M^{\vee}$, $m\in M$.
For a finite algebra $A$, the $A$-module $A^{\vee}$ is the dualising module,
see for example~\cite[Proposition~21.4]{eisenbud}, and we write it as
$\omega_A$.

We say that $M$ is
\emph{self-dual} if there is an isomorphism $\varphi\colon M\to M^{\vee}$ of
$S$-modules. If this is the case, then $\varphi^{\top}\colon
M = (M^{\vee})^{\vee}\to M^{\vee}$ is another isomorphism of $S$-modules. We say that
$\varphi$ is \emph{symmetric} if $\varphi = \varphi^{\top}$.
An isomorphism $\varphi$ comes with an associated nondegenerate $\kk$-bilinear map
$\Phi\colon M\times M\to \kk$,
that satisfies the condition $\Phi(sm_1,m_2) = \Phi(m_1, sm_2)$ for all $s\in S$, $m_1, m_2\in
M$. We refer to this condition as \emph{$S$-linearity}. The map is obtained by $\Phi(m_1, m_2) = \varphi(m_1)(m_2)$. Moreover, the
map $\Phi$ is symmetric if and only if $\varphi$ is symmetric.
\begin{example}[Self-duality from a functional]\label{ex:quadricsOnGorenstein}
    Let $A$ be a zero-dimensional $\kk$-algebra. The $A$-module $A$
    is self-dual if and only if $A$ is Gorenstein. Moreover, every isomorphism
    $\varphi\colon A\to \omega_A$ is symmetric. Indeed, let $\alpha :=
    \varphi(1)\in A^{\vee}$, then
    \[
        \Phi(a_1, a_2) = \varphi(a_1)(a_2) = \left( \varphi(a_1\cdot 1)
        \right)(a_2) =  \left(a_1\cdot \varphi (1)\right)(a_2) =
        (a_1\cdot \alpha)(a_2) = \alpha(a_1a_2),
    \]
    which is symmetric since $A$ is commutative. We see that $\varphi$ is
    completely determined by $\alpha$. Conversely, a functional $\alpha$ yields
    an isomorphism $\varphi$ if and only if $\alpha$ is nonzero on a socle
    element of every local algebra of $A$. The isomorphism $\varphi^{-1}\colon \omega_A\to A$
    is an element of $\Sym^2 A$.
\end{example}

\subsection{Iarrobino's symmetric decomposition}\label{ssec:symmDec}

    Iarrobino's symmetric decomposition of the Hilbert function of a
    Gorenstein algebra appears in
    particular in~\cite{iarrobino_associated_graded}. It plays the role of
    (refined) Hilbert function in the Gorenstein setup and it is widely
    used~\cite{cjn13,
    Iarrobino_Macias_Marques}. We recall it below.

    Let $(A, \mm, \kk)$ be a local Gorenstein $\kk$-algebra and fix a
    symmetric nondegenerate $A$-linear pairing $\Phi\colon A\times A\to \kk$.
    For every $i$, $j$, the pairing descends to a nondegenerate pairings
    \begin{align}\notag
        &\frac{\mm^i}{\mm^{i+1}} \times \frac{\annmm{i+1}}{\annmm{i}} \to \kk\\
        &\frac{\mm^i\cap \annmm{j+1}}{\mm^{i+1}\cap \annmm{j+1} + \mm^i \cap
        \annmm{j}} \times \frac{\mm^j\cap \annmm{i+1}}{\mm^{j+1}\cap
        \annmm{i+1} + \mm^j \cap \annmm{i}}\to \kk\label{eq:duality}
    \end{align}
    Let $s$ be the largest integer such that $\mm^s \neq 0$, it is called the
    \emph{socle degree} of $A$.
    Following~\cite[\S1]{iarrobino_associated_graded}, we denote for every
    $\delta=0,1, \ldots ,s-2$ the subquotient
    \[
        (\QHilb{\delta})_i := \frac{\mm^i\cap \annmm{s-i-\delta+1}}{\mm^{i+1}\cap \annmm{s-i-\delta+1} + \mm^i \cap
        \annmm{s-i-\delta}}
    \]
    and $\QHilb{\delta} := \bigoplus_{i} (\QHilb{\delta})_i$. For every $\delta$, the linear
    space $\QHilb{\delta}$
    is actually a graded $\gr(A)$-module and~\eqref{eq:duality} yields a canonical
    isomorphism $\QHilb{\delta}\to \QHilb{\delta}^{\vee}[-(s-\delta)]$ of graded $\gr(A)$-modules, where
    $[-(s-\delta)]$ denotes the shift in degree by $-(s-\delta)$. Therefore, the graded
    module
    \[
        \mathcal{Q}(A) := \bigoplus_{\delta=0}^{s-2} \QHilb{\delta}
    \]
    is self-dual (but not graded self-dual).

    Let $\DeltaHilb{\delta}$ be the Hilbert function of $\QHilb{\delta}$. Applying Hilbert
    function to the duality above,
    we obtain $\DeltaHilb{\delta}(i) =
    \DeltaHilb{\delta}(s-i-\delta)$ for every $0\leq i\leq s-\delta$ and $\DeltaHilb{\delta}(i) = 0$ for
    $i$ outside this range. Moreover, we have
    \[
        H_A = \sum_{\delta=0}^{s-2} \DeltaHilb{\delta}.
    \]
    This is Iarrobino's \emph{symmetric decomposition of the Hilbert
    function}.

    \begin{example}\label{ex:associatedGradedSymmetric}
        Let $A = \kk[x, y]/(xy, x^2 - y^3)$. The associated graded is $\gr(A)
        = \kk[x, y]/(xy, x^2, y^4)$ and so $H_A = (1,2,1,1)$ and $s = 3$. The
        algebra $\gr(A)$ is not Gorenstein.

        The element $x\in A$ lies in $\mm \cap (0:\mm^2)$ and yields a nonzero
        class in $(\QHilb{1})_1$, in fact $\QHilb{1} = \spann{x}$ is spanned
        by this class. The space $\QHilb{0}$ is $\spann{1, y, y^2, y^3}$, so
        that as $\gr(A)$-modules, we have $\QHilb{0}  \simeq  \frac{\gr(A)}{(x)}$ and
        $\QHilb{1}  \simeq  \frac{\gr(A)}{(x, y)}$.
        For quotients of $\kk[x, y]$, it is true in general that
        $\QHilb{\bullet}$ are principal
        $\gr(A)$-modules~\cite[\S2]{iarrobino_associated_graded}. This becomes
        false in more variables.
    \end{example}
    Compared to the usual theory of Hilbert functions, much less is known about
    the symmetric decompositions and Hilbert functions of Gorenstein algebras.
    For example, the
    Hilbert functions of Gorenstein quotients of $\kk[x, y, z]$ are not
    classified and neither are their symmetric decompositions.
    See~\cite{Jelisiejew_Masuti_Rossi} for a very partial result and~\cite{Iarrobino_Macias_Marques} for more discussion.

    \subsubsection{The case of modules}
    \newcommand{\annmmM}[1]{(0:_{\scriptscriptstyle M}\mm^{#1})}

        Self-dual modules received much less attention than Gorenstein
        algebras. Kunte investigated them in~\cite{Kunte__Gorenstein_modules_of_finite_length}.
        Iarrobino's symmetric decomposition has been generalised to self-dual
        modules~\cite{Wojtala}. The argument is smooth once one chooses the
        correct generalisation of $\annmm{i}$.

        Let $M$ be a zero-dimensional self-dual module over a local $\kk$-algebra $(A, \mm,
        \kk)$.
        For an ideal $I\subseteq A$ we define $(0:_{\scriptscriptstyle M} I)$ to be the set
        of elements of $M$ annihilated by $I$.
        Let $s$ be the largest integer such that $\mm^{s} M \neq 0$.
        We define
        \[
            (\QHilb{\delta})_i := \frac{\mm^iM\cap \annmmM{s-i-\delta+1}}{\mm^{i+1}M\cap
            \annmmM{s-i-\delta+1} + \mm^iM \cap
            \annmmM{s-i-\delta}}.
        \]
        Then,
        as proved in~\cite[Theorem~1.1]{Wojtala}, the graded $\gr(A)$-module
        $\QHilb{\delta}$ is self-dual:
        \[
            \QHilb{\delta} \simeq  \QHilb{\delta}^{\vee}[-(s-\delta)]
        \]
        and all other general properties from the algebra case generalise to
        this case. Again, we will see in~\S\ref{ssec:IarrobinoCurves} that the
        passage to $\QHilb{\bullet}$ is realised by a $\Gmult$-limit on the
        Quot analogue of the Iarrobino scheme.

\section{Completed quadrics}\label{sec:completedQuadrics}

\newcommand{\cV}{\mathcal{V}}%
\newcommand{\cB}{B}%
    In this chapter we discuss the variety of completed quadrics, which plays
    for the Iarrobino scheme a role similar to that played by the Grassmannian
    for the Quot scheme. It is a classical object, see~\cite{Tyrrell__complete_quadrics, Laksov__completed_quadrics,
    Concini_Procesi__Complete_symmetric_varieties, Kleiman_Thorup, DeConcini_Springer, Thaddeus__complete_collineations},
    however currently seems to be known mostly to experts.
    It is also a prominent example of a wonderful compactification.

    Despite the above sources, we were unable to find references for a number of foundational facts.  Moreover, some ideas, such as
    compatible and anticompatible bundles seem to be new.
    Hence, we give a fairly thorough discussion.
    We work in the generality of vector bundles, rather
    than vector spaces, so officially our source of reference
    is~\cite{Kleiman_Thorup} and, for~\S\ref{ssec:Tyrrell},
    also~\cite{Tyrrell__complete_quadrics}.

    Let $\cV$ be a locally free sheaf of rank $d$ over a scheme $\cB$. The variety of completed quadrics on $\cV$
    is a smooth projective morphism
    \[
        \CQ(\cV)\to \cB,
    \]
    which we recall more precisely below.
    The variety $\CQ(\cV)$ commutes with base change (by construction,
    see~\S\ref{ssec:CQconstruction} below, see also~\cite{Kleiman_Thorup}): for any
    morphism $\varphi\colon
    B'\to B$, we have a cartesian diagram
    \[
        \begin{tikzcd}
            \CQ(\varphi^*\cV) \ar[r]\ar[d]\ar[rd, "\usebox\pullback", very near
            start, phantom] &
            \CQ(\cV)\ar[d] \\
            B'\ar[r, "\varphi"] & B
        \end{tikzcd}
    \]

    \subsection{Points of $\CQ(\cV)$}\label{ssec:CQpoints}

        In this section we describe the $\kk$-points of $\CQ(\cV)$. For a nonzero quadric $q\in \Sym^2(V)$ on some
        vector space $V$, we denote by $[q]\in \mathbb{P}\Sym^2(V)$ its class
        up to rescaling. We identify such a $q$ with a
        $\kk$-linear map $q\colon V^{\vee}\to V$ that satisfies $q^{\top} =
        q$ and denote by $\coker(q)$ the space $V/\im(q)$.

        A $\kk$-point of $\CQ(\cV)$ is just a $\kk$-point of
        $\CQ(\cV|_b)$ for some $b\in \cB(\kk)$, hence we assume $B =
        \Spec(\kk)$ and write $V$ instead of $\cV$.
        We have the following description of $\kk$-points
        \begin{align}\label{eq:kpointsCQ}
            \CQ(V)(\kk) = \Big\{ ([q_{0}], [q_1], [q_2], \ldots)\ \Big|\ & q_0\colon V^{\vee} \to V,\ q_0^{\top} = q_0,\\
            & q_1\colon \ker q_0 \to \coker q_0,\ q_1^{\top} = q_1,\notag\\
            & q_2\colon \ker q_1 \to \coker q_1,\ q_2^{\top} = q_2,\notag\\
            & \phantom{q_2\colon \ker q_1 \to}\ldots\phantom{\coker q_1 \ q_2^{\top} = q_2}
        \Big\}\notag
        \end{align}
        where every $q_0$, $q_1$, $q_2$, \ldots is nonzero and the sequence
        concludes (after at most $d = \dim_{\kk} V$ steps) with a full rank quadric. We will
        abbreviate the sequence $([q_0], [q_1],\ldots)$ to $[q_{\bullet}]$.
        For indexing convenience we assume $q_{-1} = 0$, so that $\ker
        q_{-1} = V^{\vee}$ and $\coker q_{-1} = V$.

        \begin{definition}
            A $\kk$-point $[q_{\bullet}]$ is called a \emph{broken quadric}.
            For $[q_{\bullet}]\in \CQ(V)$, let $F_i :=
            \ker q_{i-1}$, so that $F_0 = V^{\vee}$. The sequence
            $F_0 = V^{\vee} \supsetneq F_1 \supsetneq F_2 \supsetneq \ldots$
            is the \emph{flag associated to $[q_{\bullet}]$}.
        \end{definition}
        For a point $[q_{\bullet}]$ with flag $F_{\bullet}$, for
        each $i\geq 0$ the quadric $q_i$ induces a \emph{full rank} quadric on
        $F_{i}/F_{i+1}$.

        The $\kk$-points above lie on the spectrum between the following two
        extreme cases presented below.
        \begin{example}[main stratum]\label{ex:mainstratum}
            If $q_0$ has full rank, then there are no other quadrics and the
            flag is trivial. In this way, the variety $(\fullrk)(\kk)$ of full
            rank quadrics up to rescaling is contained in $\CQ(V)(\kk)$. We will see
            in~\S\ref{ssec:CQconstruction} that this inclusion comes from an open embedding.
        \end{example}

        \begin{example}[most degenerate stratum]\label{ex:degstratum}
            Suppose now that $q_0$, $q_1$, \ldots all have rank
            one, so that the last quadric is $q_{d-1}$. In this case, the
            classes $[q_0]$, $[q_1]$, \ldots  are determined by the flag.
            For example if $F_1 := \ker q_0$, then the exact sequence
            \[
                \begin{tikzcd}
                    0 \ar[r] & F_1 \ar[r] & V^{\vee} \ar[r, "q_0"] & \im
                    q_0\ar[r] &
                    0
                \end{tikzcd}
            \]
            dualises to
            \[
                \begin{tikzcd}
                    0 \ar[r] & (\im q_0)^{\vee} \ar[r, "q_0^{\top} = q_0"] & V
                    \ar[r] & F_1^{\vee} \ar[r] & 0,
                \end{tikzcd}
            \]
            so $\coker q_0 = F_1^{\vee}$ and so $q_0$ is the
            unique (up to rescaling) $\kk$-linear map $q_0\colon V^{\vee}\to
            V$ with kernel $F_1$ and cokernel $F_1^{\vee}$. Let $F_{i+1} :=
            \ker q_i$ for $i=0,1, \ldots $, then the associated flag
            \[
                V^{\vee} \supseteq F_1 \supseteq F_2\supseteq \ldots \supseteq
                F_{d-1} \supseteq 0
            \]
            is a full flag on $V$. It will follow from the construction that in this way
            the flag variety becomes a closed subvariety of $\CQ(V)$.
        \end{example}

        We now introduce the notion of compatibility, which will play a role in the definition of the Iarrobino scheme.
        We keep the notation that $V$
        is a rank $d$ vector space over $\kk$.
        \begin{definition}[compatibility]\label{ref:compatibility:defNew}
            Let $x\in \End(V)$ and $[q_{\bullet}]\in \CQ(V)$. Let
            $F_{\bullet} := \ker
            q_{\bullet-1}$ be the flag associated to $[q_{\bullet}]$ and let
            $x^{\top}\in \End(V^{\vee})$ be the dual endomorphism. We say that $x$
            and $[q_{\bullet}]$ are \emph{compatible} if both conditions below
            hold
            \begin{enumerate}
                \item the endomorphism $x$ preserves the flag $F_{\bullet}$, that
                    is, $x^{\top}(F_i)\subseteq F_i$ for every $i\geq 0$.
                \item for every $i\geq0$, the induced endomorphism $\bar{x} :=
                    x^{\top}|_{F_{i}/F_{i+1}}\in \End(F_{i}/F_{i+1})$ is
                    \emph{symmetric} with respect to the full rank quadric
                    $\bar{q} := q_i|_{F_{i}/F_{i+1}}$,
                    that is
                    \[
                        \bar{q}\left( \bar{x}(-), - \right) = \bar{q}(-,
                        \bar{x}(-)).
                    \]
            \end{enumerate}
            Compatible endomorphisms for a given $[q_{\bullet}]$ form a $\kk$-vector
            subspace of $\End(V)$.
        \end{definition}
        \begin{definition}\label{ref:anticompatibility:defNew}
            In the notation of
            Definition~\ref{ref:compatibility:defNew}, we say that $x$ and
            $[q_{\bullet}]$ are \emph{anticompatible}, if the conditions from
            the definition hold with symmetric replaced by
            \emph{antisymmetric}.
        \end{definition}
        We define the $\kk$-points of the (anti)compacible bundle $\comp\to
        \CQ(V)$, $\anticomp\to \CQ(V)$ to be the following subsets of $\End(V)
        \times \CQ(V)$:
        \begin{align}\label{eq:kpointsComp}
            \comp(\kk) &= \left\{ (x, [q_{\bullet}])\ |\ \mbox{compatible}
        \right\}\\\label{eq:kpointsantiComp}
            \anticomp(\kk) &= \left\{ (x, [q_{\bullet}])\ |\ \mbox{anticompatible} \right\}
        \end{align}
        \begin{example}[main stratum, continued]
            If $q_0$ has full rank, as in Example~\ref{ex:mainstratum}, then
            $x$ is compatible with $[q_0]$ if and only if $x$ is symmetric
            with respect to it. Same holds for anticompatibility. Hence
            $\comp(\kk)|_{[q_0]}$ is a vector space of dimension
            $\binom{d+1}{2}$, while $\anticomp(\kk)|_{[q_0]}$ has dimension
            $\binom{d}{2}$.
        \end{example}

        \begin{example}[most degenerate stratum, continued]
            When $q_0$, $q_1$, \ldots have rank one, as in
            Example~\ref{ex:degstratum}, an endomorphism $x$ is compatible with
            $[q_{\bullet}]$ if and only if it preserves the flag. Therefore,
            the fibre of $\comp(\kk)$ over $[q_{\bullet}]$ is the Borel Lie
            subalgebra in $\End(V)$. In down to earth terms, up to a
            choice of basis, the fibre identifies with upper triangular matrices. In
            particular, the dimension is $\binom{d+1}{2}$.

            For an anticompatible $x$ however, the antisymmetry
            condition is not vacuous and such endomorphisms form the radical of
            the Borel Lie subalgebra, or, in down to earth terms, the
            \emph{strictly} upper triangular matrices. The dimension is
            $\binom{d}{2}$.
        \end{example}
        \begin{example}\label{ex:compatibilityVariation}
            Let us give an example where $[q_{\bullet}]$ varies. Let $d = 2$ and
            $V = \spann{e_1, e_2}$.
            Take the family $[q_{\bullet}(\lambda)]$ indexed by $\lambda\in
            \kk$, given by $q_0 =
            e_0^2+\lambda e_1^2$ for $\lambda\in \kk$. For $\lambda\neq 0$
            this is full rank, while at $\lambda = 0$ we need to put $q_1 =
            e_1^2$, which is unique up to scalars. Compatible
            operators $\comp|_{[q_{\bullet}(\lambda)]}$ yield a
            one-parameter family of subspaces
            \[
                \begin{pmatrix}
                    a_{11} & a_{12}\\
                    \lambda a_{12} & a_{22}
                \end{pmatrix} \subseteq \End(V).
            \]
            From the above description of the (would be) fibers of the compatible bundle, we conclude
            that this vector bundle (once constructed) provides
            degenerations of the space of symmetric matrices to the upper
            triangular ones. The construction itself is given in
            Proposition~\ref{ref:compatible:prop} below.
        \end{example}

        Now we describe the fibres $\comp|_{[q_{\bullet}]}$,
        $\anticomp|_{[q_{\bullet}]}$ over an arbitrary $\kkbar$-point $q_{\bullet}$.
        Take a basis $(e_1, \ldots ,e_d)$ of $V$ such that
        \begin{align}\label{eq:diagonalQuadric}
            q_0 &= e_1^2 +  \ldots + e_{r_0}^2\\\notag
            q_1 &= e_{r_0+1}^2 +  \ldots + e_{r_0+r_1}^2\\\notag
            q_2 &= e_{r_0+r_1+1}^2 +  \ldots + e_{r_0 + r_1 + r_2}^2\\\notag
             \ldots\\\notag
             q_k &= e_{r_0+r_1+ \ldots + r_{k-1}+1}^2 +  \ldots + e_{d}^2
        \end{align}
        \begin{proposition}\label{ref:orthogonalComplement:prop}
            For $[q_{\bullet}]$ given in~\eqref{eq:diagonalQuadric},
            the fibres $\comp|_{[q_{\bullet}]}$ and
            $\anticomp|_{[q_{\bullet}]}$ are given by the block
            matrices
            \[
                \begin{pmatrix}
                    A_0 & * &  \ldots & * & * & *\\
                    0 & A_1 & * &  \ldots & * & *\\
                    0 & 0 & A_2 & * & \ldots & *\\
                     &&\ldots\\
                     0 & 0 & 0 & \ldots & 0 & A_{k}
                \end{pmatrix}
            \]
            where $*$ denote arbitrary matrices, while $A_0, \ldots ,A_k$ are
            symmetric (for the compatible fibre) or antisymmetric (for the
            anticompactible fibre), respectively.
            In particular, the fibres $\comp|_{[q_{\bullet}]}$ and
            $\anticomp|_{[q_{\bullet}]}$ are orthogonal complements with
            respect to the trace pairing $(x, y)\mapsto \tr(xy)$ on $\End(V)$.
        \end{proposition}
        \begin{proof}
            The description is straightforward from the definition. Consider
            now a matrix $x\in \comp|_{[q_{\bullet}]}$ that has diagonal
            blocks $A_0, \ldots ,A_k$ and $y\in \anticomp|_{q_{\bullet}}$,
            that has diagonal blocks $A_0', \ldots ,A_k'$. To prove that
            $\tr(xy) = 0$, it is enough to prove that $\tr(A_iA_i') = 0$ for
            every $i=0,1, \ldots ,k$.

            We have $(A_iA_i')^{\top} = -A_i'A_i$, so that
            \[
                \tr(A_iA_i') = \tr((A_iA_i')^{\top}) = \tr(-A_i'A_i) =
                -\tr(A_iA_i'),
            \]
            and so the trace is zero (since $1/2\in \kk$). It remains to prove that
            \[
                \dim \comp|_{[q_{\bullet}]} + \dim \anticomp|_{[q_{\bullet}]}
                = d^2,
            \]
            but this is immediate: in the diagonal blocks, the symmetric and
            antisymmetric matrices have complementary dimensions, while each
            off-diagonal block appears twice: above and below the diagonal.
        \end{proof}

        \subsection{Construction of $\CQ(\cV)$}\label{ssec:CQconstruction}
\newcommand{\fullrkbundle}{\fullrkarg{\cV}}

            We use the convention that the vector bundle associated to $\cV$
            is $\Spec_B \Sym \cV^{\vee}$, so that the global sections of this vector bundle
            and of sheaf $\cV$ are the same. By abuse of
            notation, we will use $\cV$ to denote the vector bundle as well.
            Via standard constructions, we may in particular form the
            projective bundle $\mathbb{P}\Sym^2 \cV^{\vee}$ of quadrics.

            There are few equivalent constructions of $\CQ(\cV)$. In each of them, it is rather
            easy to realise $\CQ(\cV)$ as a set, while obtaining, e.g.,  a coordinate-free description
            of the tangent space is much harder.
            We recommend~\cite{Kleiman_Thorup,
            Thaddeus__complete_collineations} for a thorough discussion.
            For us it is most convenient to use the explicit
            \emph{description as a closed subscheme}. We summarise it below.
            For any section $q\in \Sym^2 \cV$
            let $\Lambda^{\bullet} q$ denote the tuple
            \[
                \left( q, \Lambda^2 q,  \ldots , \Lambda^{d}q \right),
            \]
            where $\Lambda^i q\in \Sym^2 \Lambda^i \cV$ is the $i$-th wedge
            power of $q$ viewed as a symmetric linear map $q\colon \cV^{\vee}\to \cV$.
            Explicitly, for any local sections $v_1^*, \ldots ,v_i^*\in \cV^{\vee}$, we have
            $(\Lambda^i q)(v_1^* \wedge \ldots \wedge v_i^*)^2 =
            \det[q(v_j^*v_k^*)]_{1\leq j,k\leq i}$.

            Consider now the product
            \begin{equation}\label{eq:product}
                \mathcal{P} := \mathbb{P}(\Sym^2 \cV) \times_B \mathbb{P}\left(
                \Sym^2\Lambda^2 \cV \right)\times_B
                \ldots \times_B \mathbb{P}\left( \Sym^{2}\Lambda^{d-1}\cV
                \right) \times_B
                \mathbb{P}\left(\Sym^2\Lambda^{d} \cV\right).
            \end{equation}
            of projective spaces, where the last factor is canonically isomorphic to $B$. Let
            $\fullrkbundle$ inside $\mathbb{P}(\Sym^2\cV)$ be the open subset of
            quadrics of full rank (that is, rank $d$). The map $q\mapsto \Lambda^{\bullet}q$ gives an
            embedding of $\fullrkbundle$ into $\mathcal{P}$. (The points in the image are called
            \emph{exterior}~\cite[p.258]{Kleiman_Thorup}.) The main theorem,
            which may serve as a definition, is the following.
            \begin{proposition}[{\cite{Laksov__completed_quadrics,
                Kleiman_Thorup}}]
                The variety of completed quadrics $\CQ(\cV)$ is the closure of
                the variety
                $\fullrkbundle$ in $\mathcal{P}$.
            \end{proposition}
            From the description, it is not clear why $\CQ(\cV)$ is smooth. On
            the positive side, the variety is projective by construction, so
            we may compute some examples of limits.

            \begin{example}\label{ex:limitInCQ}
                Let $V$ be a $\kk$-vector space of dimension $4$ with basis $e_1, e_2, e_3, e_4$.
                Consider a family of quadrics $q_{\lambda} := e_1^2 + \lambda e_2^2 +
                \lambda^2(e_3^2 + e_4^2)$, for $\lambda\in \kk$.

                For $\lambda\neq 0$ the quadric $q_{\lambda}$ has full rank,
                so it yields a point
                $\CQ(V)$. Explicitly, this point is the triple
                \begin{align*}
                    &\left[e_1^2 + \lambda e_2^2 + \lambda^2(e_3^2 +
                    e_4^2)\right],\\
                    &\left[\lambda (e_1\wedge e_2)^2 + \lambda^2\left( (e_1
                    \wedge e_3)^2 + (e_1 \wedge e_4)^2 \right) + \lambda^3\left(
                    (e_2\wedge e_3)^2 + (e_2\wedge e_4)^2
                    \right) + \lambda^4 (e_3\wedge e_4)^2\right],\\
                    &\left[\lambda^3\left((e_1\wedge e_2\wedge e_3)^2+(e_1\wedge e_2\wedge
                    e_4)^2\right) + \lambda^4\left((e_1 \wedge e_3\wedge
                    e_4)^2\right) + \lambda^5(e_2\wedge e_3\wedge e_4)^2\right]
                \end{align*}
                Let us take the limit with $\lambda\to 0$ in $\mathcal{P}$.
                It is computed coordinate-wise. In our case the
                limit is the triple
                \[
                    \left[e_1^2\right],\quad \left[(e_1\wedge e_2)^2\right],\quad \left[(e_1\wedge e_2\wedge e_3)^2 +
                    (e_1\wedge e_2\wedge e_4)^2\right].
                \]
                The first factor is the ``usual'' limit in $\mathbb{P}\Sym^2 V$, but the
                remaining ones contain much more information.
            \end{example}

            \subsubsection{Points of completed quadrics in the exterior power
            construction}\label{ssec:pointsInExterior}

            For completeness, we compare
            the description of $\kk$-points
            given in~\S\ref{ssec:CQpoints} and the one coming from the
            construction.
            Again, we reduce from $\CQ(\cV)$ to
            $\CQ(V)$, where $V$ is a $d$-dimensional $\kk$-vector space.

            For a broken quadric $[q_{\bullet}]$ we define a $\kk$-point of
            $\mathcal{P}$ defined in~\eqref{eq:product}, as the tuple
            \begin{equation}\label{eq:point}
                [q_0], [\Lambda^2 q_0],  \ldots , [\Lambda^{\rk q_0-1} q_0],\
                [\det q_0],\
                [\det q_0
                \wedge q_1],\ [\det q_0 \wedge \Lambda^2 q_1],\
                \ldots , [\det q_0 \wedge \det q_1
                \wedge q_2], \ldots
            \end{equation}
            where $\det q_{i} := \Lambda^{\rk q_i} q_i$ for every $i\geq 0$.
            The tuple~\eqref{eq:point} lies in $\CQ(V)$ and every $\kk$-point of $\CQ(V)$ has this form.
            We omit the details due to space considerations and direct the reader to the arXiv version of this article.

            \begin{example}\label{ex:limitInCQcd}
                In the setup of Example~\ref{ex:limitInCQ}, the limiting point corresponds to
                the broken quadric with a flag
                \[
                    V^{\vee} \supseteq F_1 = \spann{e_2^*, e_3^*, e_4^*}\supseteq F_2 =
                    \spann{ e_3^*,
                    e_4^*}\supseteq \{0\}\]
                    and quadrics $q_0 = e_1^2$, $q_1 = e_2^2$, $q_2 = e_3^2 +
                    e_4^2$.
            \end{example}

            \begin{example}
                Let $q_0 + t^\nu q_1\in \Sym^2(V)[\![t]\!]$ be a quadric on $\kkt$, which is
                generically of full rank.
                Section~\ref{ssec:pointsInExterior}
                might suggest that every $q_0 +
                t^{\nu} q_1 + t^{\nu+1}( \ldots )$ has the same limit as $q_0
                + t^{\nu}
                q_1$, but this is wrong.
                Indeed, take $V = \spann{e_1, e_2, e_3, e_4}$ four-dimensional and
                $q(t) := e_1^2 + e_2^2 + t^2(e_1e_3+e_2e_4)\in \Sym^2 V[\![t]\!]$. This
                has limit $([e_1^2+e_2^2],\, [e_3^2+e_4^2])$. In contrast, the
                limit of the quadric $q(t) + t^3e_3^2$ is $([e_1^2+e_2^2],\,
                [e_3^2],\,[e_4^2])$.
            \end{example}

\subsection{Tyrrell's affine patches}\label{ssec:Tyrrell}

\newcommand{\Tyr}[1]{\operatorname{Tyr}_{#1}}

Tyrrell~\cite{Tyrrell__complete_quadrics} produces a beautiful explicit open cover of
$\CQ(V)$ by affine spaces $\mathbb{A}^{\binom{d+1}{2}}$. The cover in practice
is much more convenient that the full construction.

Each patch in the cover depends on a chosen
full flag
\begin{equation}\label{eq:flag}
    0 \subseteq G_{d-1} \subseteq G_{d-2} \subseteq \ldots \subseteq G_1
    \subseteq V^{\vee}.
\end{equation}
We use the notation $G_{\bullet}$ to stress that this will \emph{not} be the
flag associated to a broken quadric.
We also fix a suitable basis
$e_{\bullet} := (e_1, \ldots ,e_d)$ of $V$, such that $G_{d-i} =
\spann{e_1^*, \ldots , e_i^*}$ for every $i=1,2, \ldots ,d-1$.
The coordinates on the patch will depend on the choice of basis. The patch
itself is independent of the choice of the basis, but depends on the choice of the flag.

Consider the open
subset $\Tyr{e_{\bullet}}\subseteq \CQ(V)$ which consists of the elements $(p_1, \ldots ,p_d)\in
\CQ(V)$ such that
\[
    p_i(\det G_i, \det G_i) = p_i\left(e_1^*\wedge \ldots \wedge e_i^*,
    e_1^*\wedge  \ldots \wedge e_i^*\right)
\]
is nonzero for every $i=1, \ldots ,d$. In other words, in bases, we consider
symmetric matrices with nonvanishing top left principal minors.
Fix an affine space of dimension $\binom{d+1}{2}-1$ with coordinates $y_1,
\ldots ,y_{d-1}$, $x_{12}$, $\ldots$ ,$x_{1d}$, $x_{23}$, $\ldots$, $x_{2d}$,
$\ldots$, $x_{d-1d}$. Let $U$ be the unipotent radical of the Borel group in
$\GL(V^{\vee})$ which fixes the flag~\eqref{eq:flag}. In coordinates, we have
\[
    U = \begin{pmatrix}
        1 & x_{12} & x_{13} & \ldots & x_{1d}\\
        0 & 1 & x_{23} & \ldots & x_{2d}\\
        && \ldots &&\\
        0 & 0 & 0 & 1 & x_{d-1d}\\
        0 & 0 & 0 & 0 & 1
    \end{pmatrix}
\]
that is, a generic unipotent upper-triangular matrix. Let $D$ be the diagonal matrix
with entries
\[
    1,\ y_1,\ y_1y_2,\ y_1y_2y_3, \ldots ,\ y_1 \ldots y_{d-1}.
\]
Consider the $i$-th wedge power $\Lambda^i D$.
It is diagonal and, by construction of $D$, it has all entries divisible by the entry
$y_1^{i-1}y_2^{i-2} \ldots y_{i-1}$ corresponding to $(e_1^*\wedge \ldots \wedge
e_i^*)^2$. Let $D^{(i)}\in \Sym^2 \Lambda^{i}V$ be
given by
\[
    D^{(i)} := \frac{1}{y_1^{i-1}y_2^{i-2} \ldots y_{i-1}}\cdot \Lambda^{i} D,
\]
so that the entry of $D^{(i)}$ corresponding to $(e_1^*\wedge \ldots \wedge
e_i^*)^2$ is $1$.  Finally, let
\begin{equation}\label{eq:TyrrellsGenericPoint}
    \left([U^{\top}D^{(1)}U],\ [U^{\top}D^{(2)}U],\ [U^{\top}D^{(3)}U],\  \ldots\right)\in
    \prod_{i=1}^{d-1} \mathbb{P}\Sym^2 \Lambda^i V.
\end{equation}
Thanks to the unit entry, this gives a well defined map from $\Spec\left(
\kk[y_1, \ldots ,y_{d-1}, x_{ij}\ |\ i<j ]\right)$ to the product. When $y_1$,
    $y_2$, \ldots, $y_{d-1}$ are all nonzero, the image is in the locus of full rank
    quadrics, hence the full image is contained in $\CQ(V)$.

\begin{proposition}[{Tyrrell~\cite{Tyrrell__complete_quadrics},
        see also~\cite[Lemma~(7.8)]{Kleiman_Thorup}}]\label{ref:Tyrrell:prop}
        The map
        \[
            \Spec\left( \kk[y_1, \ldots ,y_{d-1}, x_{ij}\ |\ i<j ]\right)\to \CQ(V)
            \]
            is an open immersion and its image is exactly
        $\Tyr{e_{\bullet}}$. Moreover, iterating over all the
        flags~\eqref{eq:flag} this yields
        an open cover of $\CQ(V)$.
\end{proposition}
We can identify $\Spec\left( \kk[y_1, \ldots ,y_{d-1}, x_{ij}\ |\ i<j ]\right)$
with $U \times \Spec(\kk[y_1, \ldots ,y_{d-1}])$. In this way, the open
immersion is $U$-equivariant, where $U$ is the group of strictly
upper-triangular matrices.
The closure of $\Spec(\kk[y_1, \ldots ,y_{d-1}])$ in $\CQ(V)$ is the permutohedral variety,
which appears prominently in Huh's work on Chow rings of
matroids.

Each Tyrrell's affine patch can be reconstructed as a dominant cell of the
\BBname{} decomposition~\cite{BialynickiBirula__decomposition} of $\CQ(V)$ for
a $\Gmult$-action which is diagonal in coordinates $e_1, \ldots .,e_d$ and has
suitable weights. The patches appear also as translates of the image of map $\Gamma$
in~\cite[Proposition~3.8]{DeConcini_Springer}.

\subsection{Representation theory and the construction of
(anti)compatible bundles}\label{ssec:compatibleBundle}

\newcommand{\boundary}{\partial\!\CQ(V)}%
\newcommand{\tanglog}{T_{\CQ(V),\,\boundary}}%
\newcommand{\cotanglog}{T_{\CQ(V),\,\boundary}^{\vee}}%

    Let $V$ be a vector space of dimension $d$, as usual.
    The variety $\fullrk$ is homogeneous under the $\PGL(V)$-action and for a full
    rank quadric $q$, its stabiliser $O(q) \subseteq \PGL(V)$ is equal to
    $\PGL(V)^{\sigma}$ for an involution $\sigma$ of $\PGL(V)$. Indeed, if we
    use $q$ to identify $V$ and $V^{\vee}$, so that $q$ is diagonal, then
    \[
        O(q) = \left\{ g\in \PGL(V)\ |\ g g^{\top} = 1 \right\},
    \]
    so the involution is $\sigma(g) = (g^{\top})^{-1}$. In such a setup, the
    variety $\fullrk$ admits a \emph{wonderful compactification} constructed
    by de Concini and
    Procesi~\cite{Concini_Procesi__Complete_symmetric_varieties}. A
    great survey of wonderful compactifications is~\cite{Pezzini__Lectures}. A source very well-adapted
    to our needs is a survey by B\u{a}libanu~\cite{Balibanu}.  Both surveys
    assume that $\kk$ has characteristic zero. Below, we explicitly
    construct the structures appearing in the wonderful compactification, hence
    we do not impose characteristic zero (but always assume characteristic
    $\neq 2$). Outside characteristic zero, a nice source from a more
    representation-theoretic perspective is~\cite{DeConcini_Springer}.

    Let $\boundary = \CQ(V) \setminus \fullrk$ be the \emph{boundary
    divisor}, with its reduced structure. This divisor splits into $d-1$ components $D_1, \ldots
    D_{d-1}$, where $D_i$ parameterises broken quadrics $[q_{\bullet}]\in
    \CQ(V)$ such
    that ``there is a break at $i$'', that is, there is an index $j$ such that
    \[
        \rk q_0 + \rk q_1+  \ldots + \rk q_j = i.
    \]
    In terms of Tyrrell's patches
    (\S\ref{ssec:Tyrrell}) the divisor $D_i$ is locally given by $y_i = 0$. In
    particular, $\boundary$ is Cartier and has simple normal crossings.

    It is natural to consider the logarithmic tangent sheaf $\tanglog$ on $\CQ(V)$,
    whose sections are vector fields tangent to
    $\boundary$, see~\cite[\S3]{Balibanu}. Let $\sll V = \left\{ x\in
        \End(V)\ |\ \tr(x) = 0
    \right\}$ be
    the Lie algebra of $\PGL(V)$. The vector fields $\sll V \subseteq H^0(T_{\CQ(V)})$
    preserve $\boundary$, hence we obtain a map
    \[
        \rho\colon \OO_{\CQ(V)}\otimes_{\kk} \sll V \to \tanglog.
    \]
    \begin{proposition}[anticompactible bundle via
        representations]\label{ref:anticompatible:prop}
        The map $\rho$ is surjective. Its kernel is a locally free sheaf of
        rank $\binom{d}{2}$ which corresponds to a bundle with $\kk$-points as
        described in~\eqref{eq:kpointsantiComp}. This bundle is called
        the \emph{anticompatible bundle} on $\CQ(V)$ and denoted by
        $\anticomp$, so that we have an exact sequence
        \begin{equation}\label{eq:anticomp}
            0 \to \anticomp \to \OO_{\CQ(V)}\otimes_{\kk} \sll V \to
            \tanglog\to 0.
        \end{equation}
    \end{proposition}
    \begin{proof}
        We will check surjectivity on Tyrrell's patches. We keep the notation
        from \S\ref{ssec:Tyrrell}. By definition, the patch depends on a basis
        $e_{\bullet} = (e_1, \ldots ,e_d)$ of $V$ and it has the form
        $\Tyr{e_{\bullet}} = U\times \mathbb{A}^{d-1}$, where the coordinates on $U$ are
        $(x_{ij})_{1\leq i<j\leq d}$ and the coordinates on
        $\mathbb{A}^{d-1}$ are $y_1, \ldots ,y_{d-1}$. The restriction of the logarithmic tangent
        sheaf $(\tanglog)|_{\Tyr{e_{\bullet}}}$ is free, generated by
        $(\partial_{x_{ij}})_{1\leq i < j\leq d}$ and
        $(y_i\partial_{y_i})_{1\leq i\leq d-1}$. The free $U$-action on
        $\Tyr{e_{\bullet}}$ yields an injection $\mathfrak{u}\into
        H^0(\Tyr{e_{\bullet}},
        \tanglog)$, whose image is the span
        \[
            \spann{\partial_{x_{ij}}\ |\ 1\leq i < j\leq d}.
        \]
        The $U$-action is part of the $\PGL(V)$-action, so
        the $\mathfrak{u}\otimes_{\kk} \OO_{\Tyr{e_{\bullet}}}$
        part of the logarithmic tangent sheaf is in the image of $\rho$. Consider the Lie algebra $\mathfrak{t}$ of the diagonal torus
        in $\PGL(V)$. Fix the basis of $\mathfrak{t}$ given by
        \[
            t_i = (1, 0,  \ldots 0, -1, 0, \ldots 0)
        \]
        for $i=1, \ldots ,d-1$.
        To prove surjectivity of $\rho$,
        we can work modulo $\mathfrak{u}$, so it is enough to understand the
        action of $\mathfrak{t}$ on the universal quadric
        \begin{equation}\label{eq:univOfAd}
            (e_1)^2 + y_1(e_2)^2 + y_1y_2(e_3)^2 + \ldots + y_{1} \ldots
            y_{d-1}e_d^2
        \end{equation}
        on $\mathbb{A}^{d-1}$.
        The element $(\lambda,1, \ldots, 1, \lambda^{-1},1, \ldots ,1)$ in the
        diagonal torus acts on~\eqref{eq:univOfAd} by multiplying the entry
        $y_1 \cdot  \ldots \cdot  y_i\cdot (e_i)^2$ by
        $\lambda$ and leaving all other entries intact.
        This implies that we have
        \[
            \rho(t_i) = \begin{cases}
                y_{i}\partial_{y_i} - y_{i+1}\partial_{i+1} \mod
                \rho(\mathfrak{u}) & \mbox{ for }
                i <
                d\\
                y_{i}\partial_{y_i}\mod
                \rho(\mathfrak{u})\mbox{ for } i = d
            \end{cases}
        \]
        and surjectivity of $\rho$ follows. Thus the kernel $\ker \rho$ is
        locally free of rank $d^2 - 1 - (\binom{d+1}{2} - 1) = \binom{d}{2}$.
        We identify it with the corresponding vector bundle $\ker \rho \subseteq
        \CQ(V) \times \End(V)$.
        It remains to prove that $\kk$-points of this bundle are exactly the
        anticompatible operators from
        Definition~\ref{ref:anticompatibility:defNew}.
        The
        map $\rho$ and the notion of anticompatibility are both
        $U$-equivariant, so it is enough to prove this on
        $\mathbb{A}^{d-1} \into \Tyr{e_{\bullet}}$.
    Take variables $a_{ij}$ for $1\leq i \leq j \leq d$, which will be the
    coordinates on our bundle and consider the subbundle of
    $\mathbb{A}^{d-1} \times \End(V)$ given by the \emph{transpose} of the matrix
    \[
        \left[ \begin{array}{c c}
            -a_{ij} & \mbox{for } i>j\\
            0 & \mbox{ for } i = j\\
            y_iy_{i+1} \ldots y_{j-1}\cdot a_{ij} & \mbox{for } i< j
    \end{array}\right] =
        \begin{pmatrix}
            0 & y_1a_{12} & y_1y_2a_{13} & y_1y_2y_3a_{14}& \ldots &y_1 \ldots
            y_{d-1} a_{1d}\\
            -a_{12} & 0 & y_2a_{23} & y_2y_3a_{24} & \ldots & y_2  \ldots
            y_{d-1} a_{2d}\\
            -a_{13} & -a_{23} & 0 & y_3a_{34} & \ldots & y_3 \ldots y_{d-1} a_{3d}\\
            &&  \ldots & \ldots &&\\
            -a_{1,d-1} & -a_{2,d-1} & -a_{3,d-1} & \ldots & 0 &
            y_{d-1}a_{d-1,d}\\
            -a_{1d} & -a_{2d} & -a_{3d} & \ldots & -a_{d-1,d} & 0
        \end{pmatrix}
    \]
    Denote the matrix by $M_{\anticomp}\in \mathbb{A}^{d-1} \times
    \End(V^{\vee})$.
    Consider the quadric~\eqref{eq:univOfAd} and let $Q$ be the associated symmetric matrix.
    We have $M_{\anticomp}^{\top} Q = -Q M_{\anticomp}$, since the $(i,j)$ entry of both sides is
    \[
        \begin{cases}
            y_{1} \ldots y_{j-1} a_{ij} & \mbox{ if } i < j\\
            0 & \mbox{ if } i =j\\
            -y_{1} \ldots y_{j-1} a_{ij} & \mbox{ if } i > j.
        \end{cases}
    \]
    This identifies the $\ker \rho$ restricted to $\mathbb{A}^{d-1}$.
    To check that the subbundle above describes exactly the anticompatible
    endomorphisms, pick a broken quadric $[q_{\bullet}]$ in
    $\mathbb{A}^{d-1}$. Thanks to $M_{\anticomp}^{\top} Q = -Q M_{\anticomp}$, we
    see that the endomorphisms preserve $\ker q_0$ and yield antisymmetric matrices in $V/\ker q_0$.
    Suppose that $q_0$ has rank $r < d$. This means that $y_r = 0$.
    The matrix $M_{\anticomp}|_{y_r = 0}$ is a $2\times 2$ block matrix. The
    lower-right block is the analogue of $M_{\anticomp}$ on the space
    $\spann{e_{r+1}, \ldots ,e_d}$. By induction, it is anticompatible with
    $\left( [q_1], [q_2], \ldots  \right)$ and the claim follows.
    \end{proof}

    \begin{remark}
        In characteristic zero, for every wonderful compactification, the
        corresponding map $\rho$ is surjective by a result of
        Pezzini~\cite[Proposition~4.2]{Pezzini__Lectures}, and the kernel of
        $\rho$ is described by
        Brion~\cite[Proposition~2.1.2]{Brion__loghomogeneous}.
    \end{remark}

    \begin{proposition}[Compatible bundle]\label{ref:compatible:prop}
        Consider the dual of the exact sequence~\eqref{eq:anticomp} with an
        added trivial factor coming from $\gl V = \sll V \oplus \kk$:
        \[
            0\to \cotanglog \oplus \OO_{\CQ(V)}\to \OO_{\CQ(V)}\otimes_{\kk}
            \gl V \to \anticomp^{\vee} \to 0.
        \]
        The vector bundle corresponding to $\cotanglog \oplus \OO_{\CQ(V)}$
        has rank $\binom{d+1}{2}$ and
        $\kk$-points as described in~\eqref{eq:kpointsComp}. We call it the
        \emph{compatible bundle} and denote it by $\comp$.
    \end{proposition}
    \begin{proof}
        The kernel of the surjection $\OO_{\CQ(V)}\otimes_{\kk}
            \gl V \to \anticomp^{\vee}$ is, fibre by fibre, the orthogonal
            complement of $\anticomp$ with respect to the trace pairing. This
            is equal to the fibre of $\comp$ by
            Proposition~\ref{ref:orthogonalComplement:prop}.
    \end{proof}

    \begin{example}
        In the case $\dim V=2$, we see that $\CQ(V)$ is isomorphic to
        $\mathbb{P}^2 = \mathbb{P} \Sym^2 V$.
        Fix coordinates $V = \spann{e_1, e_2}$. For a quadric $q =
        \lambda_{11}e_1^2 + 2\lambda_{12}e_1e_2 + \lambda_{22}e_2^2$, the
        element
        \[
            \begin{pmatrix}
                \lambda_{12} & -\lambda_{11}\\
                \lambda_{22} & -\lambda_{12}
            \end{pmatrix}\in \End(V)
        \]
        is anticompatible with $q$, so that $\anticomp$ corresponds to the
        sheaf
        $\OO_{\mathbb{P}^2}(-1)$ and $\tanglog$ is $T_{\mathbb{P}^2}(-1)$. Using
        Proposition~\ref{ref:compatible:prop}, we obtain that $\comp$ comes
        from the sheaf $\OO_{\mathbb{P}^2}\oplus \Omega_{\mathbb{P}^2}(1)$.
    \end{example}

    \begin{remark}
        The construction of $\anticomp$ and $\comp$ generalise immediately to
        the relative setting of $\CQ(\cV)$: the surjectivity of
        $\rho$ is checked locally, hence we can assume that $\cV$ is trivial
        and the proof above applies.
    \end{remark}

    \subsection{Duality}\label{ssec:dualityOnCQ}

        It may seem that the theory involves a choice: one can work on
        $\CQ(V)$ or on $\CQ(V^{\vee})$. In
        this section we prove that both viewpoints are canonically isomorphic. This is important, since each
        of the two choices is sometimes preferable.

        A full rank quadric on a section of $\cV$ induces one on $\cV^{\vee}$ and vice versa,
        so the open subsets $\fullrkarg{\cV}$ and $\fullrkarg{\cV^{\vee}}$ are
        isomorphic. Explicitly, if $q\in \Sym^2 \cV$ is of full rank
        and we write it as a symmetric matrix $q\colon \cV^{\vee}\to \cV$, then
        $q^{-1}\colon \cV\to \cV^{\vee}$ is a symmetric matrix and it gives the
        corresponding element of $\fullrkarg{\cV^{\vee}}$.

        The main duality theorem (see~\cite[Introduction]{Kleiman_Thorup}) is
        that this isomorphism extends to an
        isomorphism
        \[
            \CQ(\cV)\to \CQ(\cV^{\vee}).
        \]
        The isomorphism is constructed as follows. Recall that $\cV$ has
        rank $d$. Let $\det \cV  := \Lambda^d
        \cV$, then for $i=0,1 \ldots ,d $ we have isomorphisms
        \[
            \det \cV \otimes_{\OO_B} \Lambda^{i} \cV^{\vee}  \simeq
            \Lambda^{d-i} \cV
        \]
        given by the contraction of wedges and, correspondingly, an
        isomorphism
        \[
            \mathbb{P}\Sym^2 \Lambda^{i} \cV^{\vee} \simeq \mathbb{P}\Sym^2\Lambda^{d-i} \cV.
        \]
        This induces an isomorphism of
        ambient varieties
        \[
            \mathbb{P}(\Sym^2\cV^{\vee}) \times_B \mathbb{P}\left( \Sym^2\Lambda^2
            \cV^{\vee} \right)\times_B
            \ldots \times_B
            \mathbb{P}\left(\Sym^2\Lambda^{d} \cV^{\vee}\right)\to
            \mathbb{P}(\Sym^2\cV) \times_B
            \ldots \times_B
            \mathbb{P}\left(\Sym^2\Lambda^{d} \cV\right).
        \]
        The isomorphism extends the map $q\mapsto q^{-1}$ on full rank
        quadrics and $\CQ(\cV)$, $\CQ(\cV^{\vee})$ are by definition the
        closures of this loci, so we obtain the
        following result.
        \begin{proposition}[duality, {\cite[(7.23)]{Kleiman_Thorup}}]
            The isomorphism above restricts to an isomorphism $\iota_\cV\colon
            \CQ(\cV^{\vee})\to \CQ(\cV)$ such that
            $\iota_{\cV^{\vee}}\circ \iota_\cV$ is the identity.
        \end{proposition}

        \begin{example}
            In the setup of Examples~\ref{ex:limitInCQ}, \ref{ex:limitInCQcd},
            the duality interchanges the flags
            \[
                V^{\vee} \supseteq \spann{e_2^*, e_3^*, e_4^*}\supseteq \spann{e_3^*,
                e_4^*} \supseteq 0\qquad\mbox{and}\qquad
                V \supseteq
                \spann{e_1, e_2}\supseteq \spann{e_1} \supseteq 0
            \]
            and maps the quadrics $q_0 = e_1^2$, $q_1 = e_2^2$, $q_2 = e_3^2+e_4^2$ to their
            inverses $(e_1^*)^2$, $(e_2^*)^2$,
            $(e_3^*)^2+(e_4^*)^2$, arranged in reversed order.
        \end{example}

        The functor of point of $\CQ(\cV)$ is given in Thorup-Kleiman~\cite[\S7]{Kleiman_Thorup}. We omit discussing it here.

\section{Compatible families and the Iarrobino scheme}\label{sec:Iarrobino}

\newcommand{\cQ}{\mathcal{Q}}%
\newcommand{\alphaCQ}{\alpha_{\CQ}}%

    With appropriately developed theory of $\CQ(\cV)$, the definition of
    compatible families and the Iarrobino scheme of a scheme $X$ becomes
    quite natural. We realise it as a special case of a more general construction
    of the compatible locus.

    \subsection{Points of compatible quadrics}

    Let $S$ be a $\kk$-algebra. Let $V$ be a $d$-dimensional vector space
    that is an $S$-module.  Suppose that $[q_{\bullet}]\in \CQ(V)$ is
    compatible with the action of every element of $S$. Let $F_{\bullet} =
    \ker q_{\bullet -1}$ be the flag associated to $[q_{\bullet}]$. As a
    consequence of compatibility, the subspaces and quadrics become $S$-linear,
    as follows.

    \begin{proposition}\label{ref:assGr:prop}
        The following are equivalent for a quadric $[q_{\bullet}]\in \CQ(V)$.
        \begin{enumerate}
            \item $[q_{\bullet}]$ is compatible with (every element of) the
                ring $S$.
            \item
                For every $i\geq0$, the subspace $F_{i}$ is an $S$-submodule
                of $V^{\vee}$
                and
                the quadric $q_i$ induces an isomorphism
                \begin{equation}\label{eq:tmpIso}
                    \left(\frac{F_i}{F_{i+1}}\right)^{\vee}\to \frac{F_i}{F_{i+1}}
                \end{equation}
                of $S$-modules that satisfies $q_i = q_i^{\top}$. Therefore,
                the associated graded $S$-module $\bigoplus_{i\geq 0}
                F_i/F_{i+1}$ is self-dual via $q_{\bullet}$.
        \end{enumerate}
    \end{proposition}
    \begin{proof}
        By Definition~\ref{ref:compatibility:defNew}, the spaces
        $F_{\bullet}$ are preserved by the action of any $x\in S$, hence are
        submodules. Fix $i\geq 0$. From the same definition, we see that $q_i$
        satisfies $q_i(x\cdot (-), -) = q_i(-, x\cdot (-))$ for every $x\in
        S$. This is tantamount to saying that the linear
        isomorphism~\eqref{eq:tmpIso}
        is $S$-linear.
    \end{proof}

    \begin{example}[The case {$S = \kk[x]/(x^d)$}]\label{ex:completedQuadricsOnJordan}
    This example will be useful later. Let $S =
    \kk[x]/(x^d)$ and $V = S^{\oplus 1}$ so that $V^{\vee} \simeq S^{\oplus
    1}$.
Take any $[q_{\bullet}]$ compatible with
$x$. By Proposition~\ref{ref:assGr:prop} each subspace $F_i\subseteq
\kk[x]/x^d$ is an ideal, so it has the form $(x^{\nu_i})/(x^d)$ for $\nu_i := \dim
V^{\vee}/F_i$.

A quadric $q_i$ is defined on $(x^{\nu_i})/(x^{\nu_{i+1}}) \simeq
\kk[x]/(x^{\nu_{i+1}-\nu_i})$, where $\nu_0 = d$. The algebra $A_i :=
\kk[x]/(x^{\nu_{i+1} - \nu_i})$ is
Gorenstein and $q_i$ has full rank, so, by
Example~\ref{ex:quadricsOnGorenstein}, the quadric $q_i$ on
$A_i$ has the form
\[
    q_i(a', a'') = \alpha_i(a'a'')
\]
for a functional $\alpha_i$ which is nonzero on the class of $x^{\nu_{i+1} -
\nu_i
-1}$. Since we consider $q_i$ up to scalars, we may even assume that
$\alpha(x^{\nu_{i+1} - \nu_i
-1}) = 1$. In this case we obtain a $(\nu_{i+1}-\nu_i-1)$-dimensional affine space
parameterising possible $q_i$.
This completes the set-theoretic description.
We will see in Theorem~\ref{ref:smoothnessCurve:thm} that the whole compatible
locus is irreducible. A willing reader can do this now by hand.
\end{example}

\begin{example}\label{ex:completedQuadricsOnSmooth}
    Let $S = \kk\times  \ldots \times \kk$ be the algebra which is the product of $d$ copies of the
    field $\kk$. Let $V = S^{\oplus 1}$ with basis $e_1, \ldots ,e_d$. From
    Example~\ref{ex:quadricsOnGorenstein}, we see that the full
    rank quadrics compatible with $S$ correspond to full rank diagonal
    matrices. The closure of such matrices in $\CQ(V)$ is the permutohedral
    toric variety. It will follow from
    Theorem~\ref{ref:smoothnessCurve:thm} that this closure coincides with the
    variety of all broken quadrics compatible with $S$.
\end{example}

    \subsection{Compatible quadrics scheme-theoretically, affine version}
    The compatibility of a single operator $x\in \End(V)$ with a given
    $[q_{\bullet}]\in \CQ(V)$ was given in
    Definition~\ref{ref:compatibility:defNew}. By
    Proposition~\ref{ref:compatible:prop}, it can be rephrased by saying that
    $x\in \End(V)$ lies in $\comp|_{[q_{\bullet}]}\subseteq \End(V)$. The
    latter formulation immediately generalises to more operators and vector
    bundles.

    Let $X = \Spec(S)$ be an affine $\kk$-scheme and $\cV$ be a coherent
    $\OO_{X\times B}$-module which is locally free of rank $d$ as a $\OO_B$-module.
    The $\OO_{X\times B}$-module structure induces a morphism $\alpha\colon S\times B\to \End(\cV)$ of
    vector bundles over $B$,
    where $\End(\cV)$ is the vector bundle corresponding to the endomorphisms
    of the locally free
    $\OO_B$-module $\cV$ and $S = H^0(\OO_X)$ so that $S\times B$ is a trivial
    vector bundle (usually of infinite rank) on $B$.
    The morphism pulls back to a morphism
    \[
        \alphaCQ = \alpha \times \id_{\CQ(\cV)}\colon S \times \CQ(\cV)\to
        \End(\cV) \times_B \CQ(\cV)
    \]
    of vector bundles on $\CQ(\cV)$.

    \begin{definition}[Compatibility, relative affine
        version]\label{ref:compatibleLocus:def}
        The \emph{scheme of compatible quadrics on $\cV$} is a closed
        subscheme $\CCQ(\cV) \into \CQ(\cV)$ given by the vanishing of the
        following morphism of vector bundles over $\CQ(\cV)$:
        \begin{equation}\label{eq:compatibility}
            S \times \CQ(\cV) \to \End(\cV) \times_B \CQ(\cV) \to
            \frac{\End(\cV) \times_B \CQ(\cV)}{\comp}.
        \end{equation}
    \end{definition}
    In words, we require that $\alphaCQ$ factors through the compatible bundle
    $\comp \into \End(\cV) \times_B \CQ(\cV)$.

    The right-hand-side of~\eqref{eq:compatibility} is a vector bundle of rank $\binom{d}{2}$ on
    $\CQ(\cV)$, see Proposition~\ref{ref:compatible:prop}. Locally on $\CQ(\cV)$, it is trivial, so that the right hand
    side is isomorphic to $\CQ(\cV) \times \mathbb{A}^{\binom{d}{2}}$ and the
    $(\dim_{\kk} S)\cdot \binom{d}{2}$
    coordinate functions cut out the scheme $\CCQ(\cV)$. The zero sets of these local equations
    agree on overlaps, giving the scheme structure globally. If $S$ is a
    $\kk$-algebra generated by $n$ elements, then only the multiplication by
    these elements needs to be taken into account and we reduce to $n\cdot
    \binom{d}{2}$ equations.

    \begin{definition}\label{ref:Iar:def}
        Let $\pr\colon \mathcal{Z} \to \HilbdX$ be the universal family and
        $\mathcal{U} = \pr_*\OO_{\mathcal{Z}}$ be the corresponding rank $d$
        vector bundle. The scheme
        $\CCQ(\mathcal{U})\to \HilbdX$ is called the \emph{Iarrobino scheme of $d$ points
        on $X$} and denoted $\Iar_d(X)$.
    \end{definition}

    \begin{proposition}[Points of the Iarrobino scheme]\label{ref:Iarkpoints:prop}
        The $\kk$-points of $\Iar_d(X)$ correspond bijectively to pairs
        $[Z_{\bullet}], [q_{\bullet}]$, where $Z_{\bullet}$ is a flag
        \[
            X\supseteq Z_0\supsetneq Z_1\supsetneq Z_2\supsetneq \ldots 
        \]
        of zero-dimensional subschemes of $X$ such that $\deg Z_0 = d$,
        and $[q_{\bullet}] = \left( [q_0], [q_1], \ldots  \right)$ are classes
        of $\OO_X$-module isomorphisms
        \[
            q_i\colon \left(\frac{I_{Z_{i+1}}}{I_{Z_i}}\right)^{\vee}\to
            \frac{I_{Z_{i+1}}}{I_{Z_i}},
        \]
        which satisfy $q_i^{\top} = q_i$.
    \end{proposition}
    \noindent We stress that both the length of the tuple $Z_0 \supsetneq
    Z_1\supsetneq
    \ldots $ and the degrees of $Z_1, Z_2, \ldots$ are allowed to vary.

    \begin{proof}
        We write $\OO_{Z_i}$ for $H^0(\OO_{Z_i})$ for $i=0,1 \ldots$.

        By construction of $\Iar_d(X)$, its $\kk$-points correspond
        bijectively to broken quadrics $[q_{\bullet}]$ on $\OO_{Z_0}$
        compatible with $\OO_{Z_0}$-module structure, where  $Z_0\subseteq X$
        ranges over zero-dimensional subschemes of degree $d$.

        We have a sequence of cokernels
        \[
            \OO_{Z_0}\onto \coker q_0 \onto \coker q_1 \onto \ldots,
        \]
        and we define, for $i\geq 1$, the subscheme $Z_i$ as given by the cokernel
        $\OO_{Z_0}\to \coker q_{i-1}$, so that we have
        \[
            \begin{tikzcd}
                \omega_{Z_0}\ar[d, "q_0"] & \ar[l, hook'] \omega_{Z_1}\ar[d, "q_1"] &
                \ar[l, hook'] \omega_{Z_2}\ar[d, "q_2"] &\ar[l, hook'] \ldots\\
                \OO_{Z_0} \ar[r, two heads] & \OO_{Z_1} \ar[r, two heads] &
                \OO_{Z_2} \ar[r, two heads] & \ldots
            \end{tikzcd}
        \]
        and every quadric $q_i$ yields a symmetric isomorphism
        $(I_{i+1}/I_{i})^{\vee} \simeq \omega_{Z_{i}}/\omega_{Z_{i+1}}\to I_{i+1}/I_{i}$.
    \end{proof}

    Just as the Hilbert schemes of points generalise to Quot schemes of
    points, the Iarrobino scheme generalises to its Quot counterpart.
    \newcommand{\cM}{\mathcal{M}}%
    \begin{definition}\label{ref:Quot:def}
        Let $\cE$ be a coherent sheaf on $X$ and consider the Quot scheme of
        points $\Quot_d(\cE)$ with its universal quotient sheaf
        \[
            \cE \otimes_{\kk} \OO_{\Quot_d(\cE)}\onto \cM.
        \]
        The \emph{self-dual Quot scheme of $d$ points} is $\CQuot_d(\cE)$
        given by the compatibility locus $\CCQ(\cM)\to \Quot_d(\cE)$.
    \end{definition}

    There is a Iarrobino-to-Hilbert map
    \[
        \tau_{X}\colon \Iar_d(X)\to \OHilb_d(X),
    \]
    which forgets about $[q_{\bullet}]$. The fibre $\tau_X^{-1}(Z\subseteq X)$
    is the scheme of quadrics compatible with $\OO_Z$, hence is independent of
    the embedding of $Z$ into $X$.

    Below, we investigate $\Iar_d(C)$ for a curve $C$ much more precisely.
    Before that, we extend the definition from the affine to quasi-projective
    setup.

    \subsubsection{Compatible quadrics, from affine to quasi-projective version}

    If $X$ is quasi-affine, open in an affine $\Xbar$, then $\OHilb_d(X)$ is
    open in $\OHilb_d(\Xbar)$, hence we obtain $\Iar_d(X)$ open in
    $\Iar_d(\Xbar)$.

    Let $X$ now be a quasi-projective $\kk$-scheme and $\cV$ a coherent $\OO_{X\times
    B}$-module which is locally free of rank $d$ as a $B$-module. For every
    point $b\in B$, the support of $\cV|_b$ is finite as a subset of $X\times
    b$, hence by
    quasi-projectivity, there exists a hypersurface $H\subseteq X$ such that
    $H\subseteq X$ does not intersect the image of the support of $\cV|_b$. It
    follows that for some open subset $B'\subseteq B$ the support of $\cV$ is
    disjoint from $X\times B'$ and hence $\cV$ can be viewed as a coherent
    $\OO_{X\setminus H \times B'}$-module. The scheme $X\setminus H$ is
    quasi-affine. Applying the construction from above
    to $\cV$, we obtain $\CCQ(\cV|_{B'})$, which is independent up to isomorphism of
    the choice of $H$. Gluing such families for all $b\in B$, we obtain the
    desired scheme.
    Similarly one obtains the self-dual Quot scheme for any coherent sheaf on
    a quasi-projective $X$.

    \subsubsection{Stacks}\label{ssec:stacks}

        As a general reference for stacks, we use~\cite{Olsson,
        stacks_project}.
        Let $X$ be a quasi-projective scheme. The stack $\Coh_d X$ is a fibered
        category over $\Spec(\kk)$-schemes, whose fibre over a scheme $B$ is
        the groupoid of coherent $\OO_{X\times B}$-modules $\mathcal{M}$ such that
        $\mathcal{M}$ is finite flat of degree $d$ over $B$.
        The universal family on the Quot scheme yields a map
        $\Quot_d(\OO_{X}^{\oplus d})\to \Coh_d X$ that is a smooth cover and
        so $\Coh_X$ is an Artin stack.

        \begin{definition}
            The Artin stack $\CQCoh_d X\to \Coh_d X$ is the compatible locus
            $\CCQ(\mathcal{M})$ for the universal family $\mathcal{M}$ on the stack $\Coh_d X$.
        \end{definition}
        In other
        words, for a scheme $B$, the fibre of $\CQCoh_d X$ is the groupoid
        of pairs $(M, [q_{\bullet}])$, where is a coherent $\OO_{X \times
        B}$-module, locally free of rank $d$ as an $\OO_B$-module, and
        $[q_{\bullet}]\colon B\to \CCQ(M)$ is a morphism, where $\CCQ(-)$ is
        the scheme of compatible quadrics as in
        Definition~\ref{ref:compatibleLocus:def}, so that $[q_{\bullet}]$
        corresponds to a family of compatible broken quadrics over $B$.
        The morphism $\CQCoh_d X\to \Coh_d X$ forgets about $[q_{\bullet}]$.

        Completed quadrics and compatibility commute with base change, so we have a
        cartesian diagram
        \[
            \begin{tikzcd}
                \CQuot_d(\OO_X^{\oplus d}) \ar[r] \ar[rd, "\usebox\pullback", very near
                start, phantom]\ar[d, two heads, "\mathrm{sm}"] & \Quot_d(\OO_X^{\oplus d})
                \ar[d, two heads, "\mathrm{sm}"]\\
                \CQCoh_d X \ar[r] & \Coh_d X
            \end{tikzcd}
        \]
        which proves that $\CQCoh_d X$ is indeed an Artin stack.

        We now pass to the stack of finite algebras. The Artin stack $\Alg_d$ is the fibered
        category over $\Spec(\kk)$ whose fiber over a scheme $B$ is the
        groupoid of finite flat degree $d$ maps $Y\to B$. It admits a smooth
        cover $\OHilb_d(\mathbb{A}^{d-1})\to \Alg_d$.

        The Artin stack $\CQAlg_d\to \Alg_d$ is defined by considering pairs
        $(\pi\colon Y\to B, [q_{\bullet}])$, where $[q_{\bullet}]\in \CQ(\pi_*
        \OO_Y)$ is compatible with the action of $Y$. This requires a slight
        extension of the definition of compatibility.

    \subsection{The unbroken locus and the main component of the Iarrobino
    scheme}\label{ssec:unbroken}

        In this section we assume that $X$ is irreducible.
        In this case, the Hilbert scheme has a distinguished irreducible component, the smoothable
        component $\OHilb_d^{\sm}(X)$, that parameterises tuples of points and
        more general smoothable schemes, see~\cite{Bertin__punctual_Hilbert_schemes}. We
        now describe a similar irreducible component in $\Iar_d(X)$. Recall that $\HilbdGororX \to \HilbdGorX$ is a
        smooth morphism (with $(d-1)$-dimensional fibres) and the space
        upstairs parameterises pairs $(Z\subseteq
        X, [q])$, where $q\colon \omega_Z\to \OO_Z$ is an orientation.

        \begin{definition}\label{ref:unbroken:def}
            The \emph{unbroken locus} $\Iar_d^{\Gor}(X)$ in $\Iar_d(X)$ is the locus
            where the broken quadric $[q_{\bullet}]$ is not broken, that is,
            where $q_0$ has full rank.
        \end{definition}
        \begin{lemma}\label{ref:unbrokenopen:lem}
            The unbroken locus $\Iar_d^{\Gor}(X)$ is open in $\Iar_d(X)$ and isomorphic to
            $\HilbdGororX$. If $X$ is smooth and $\dim X \leq 3$, then this
            locus is smooth.
        \end{lemma}
        \begin{proof}
            Openness follows for the openness of the locus of full rank quadrics in
            $\CQ(\cV)$ for every $\cV$, see~\S\ref{ssec:CQconstruction}.
            By
            Proposition~\ref{ref:Iarkpoints:prop}, the $\kk$-points of this locus correspond to
            pairs $(Z, [q])$, for $Z\subseteq X$ a zero-dimensional \emph{Gorenstein} subscheme of
            degree $d$ and $q\colon \omega_{Z}\to \OO_Z$ an orientation, hence
            $\Iar_d(X)\to \HilbdGorX$ is smooth and surjective. It follows
            that $\Iar_d^{\Gor}(X)$ is isomorphic to $\HilbdGororX$.

            For $\dim X\leq 3$, smoothness follows from the Buchsbaum-Eisenbud
            theorem~\cite{BuchsbaumEisenbudCodimThree}, as proved
            in~\cite{kleppe_roig_codimensionthreeGorenstein}.
        \end{proof}

    \newcommand{\Iarzero}{\left(\Iar_d(X)\right)^{\circ}}
    \begin{definition}\label{ref:smoothableComponent:def}
            The \emph{locus of tuples of points} $\Iarzero$ is the intersection of
            $\Iar_d^{\Gor}(X)$ with the locus $Z\subseteq X$ such that $Z$ is
            reduced. The locus $\Iarzero$ is irreducible, open in $\Iar_d(X)$
            and has
            dimension $(1+\dim X)d - 1$. The closure of $\Iarzero$
            is called the \emph{smoothable component} and denoted
            $\Iar^{\sm}_d(X)$.
        \end{definition}

    \subsection{Iarrobino scheme of a smooth curve is smooth}\label{ssec:curves}

    In this section, we discuss $\Iar_d(C)$ for $C$ a smooth curve. Here, $\Iar_d(C)$ is smooth and so is the symmetric
    product $\OHilb_d(C)$, yet the natural map
    $\Iar_d(C)\to \OHilb_d(C)$ is not smooth.

    Similarly to the Hilbert scheme, the case of affine line $C =
    \mathbb{A}^1$ is concrete and serves as a local model.
    Recall~\cite{fogarty} that
    $\OHilb_d(\mathbb{A}^1)$ is isomorphic to $\mathbb{A}^d$, the isomorphism
    sends a finite subscheme $V(x^d + a_{d-1}x^{d-1}+ \ldots +a_0)$ to
    $(a_{d-1}, \ldots ,a_0)$. The Iarrobino scheme $\Iar_d(\mathbb{A}^1)$
    is fibred over
    $\OHilb_d(\mathbb{A}^1)$ in permutohedral varieties
    (Example~\ref{ex:completedQuadricsOnSmooth}) and their degenerations.

    It is convenient to prove smoothness of $\Iar_d(C)$ in the generality of the Quot
    schemes.
    \begin{theorem}\label{ref:smoothnessCQuotCurve:thm}
        Let $C$ be a smooth irreducible curve and $\cE$ be a locally free sheaf on $C$.
        Then the scheme $\CQuot_d(\cE)$ is smooth connected of dimension
        $d\cdot (1+\rk \cE) - 1$.
    \end{theorem}
    \begin{proof}
        Let us begin with the proof of smoothness.
        After passing to $\kk$ algebraically closed, it is enough to prove
        smoothness on $\kk$-points. Take a $\kk$-point $(M, [q_{\bullet}])$ of
        $\CQuot_d(\cE)$, where $p_C\colon \cE\onto M$ is a surjection of
        $\OO_C$-modules. Let $\{c_1, \ldots ,c_e\} \subseteq C(\kk)$ be the
        support of $M$. Passing to an open neighbourhood $U$ of the support, we
        may assume that $\cE|_{U} = \OO_{U}^{\oplus r}$ is free.

        Fix any pairwise different $\kk$-points $p_1, \ldots ,p_e\in
        \mathbb{A}^1$. By Cohen's structure theorem, there is an isomorphism
        of products of complete local rings
        \begin{equation}\label{eq:local}
            \iota\colon \prod_{i=1}^e \OOhat_{\mathbb{A}^1, p_i}\to
            \prod_{i=1}^e \OOhat_{C, c_i}.
        \end{equation}
        Let $\cE' := \OO_{\mathbb{A}^1}^{\oplus r}$.
        The $\OO_C$-module structure of $M$ is completely
        determined by the action of complete local rings, hence there exists a
        surjection $p_{\mathbb{A}^1}\colon \cE'\onto M'$ of
        $\OO_{\mathbb{A}^1}$-modules, such that $\iota^* p_{C} =
        p_{\mathbb{A}^1}$, in particular, $\iota^* M = M'$.

        The same argument applies to show that the deformations functors of
        $p_{\mathbb{A}^1}$ and $p_C$ over Artin rings are isomorphic. In
        particular, the complete local rings at $p_{C} = [\cE\onto M]\in
        \Quot_d(\cE)$ and $p_{\mathbb{A}^1} = [\cE'\onto M']\in
        \Quot_d(\cE')$ are isomorphic in a way compatible with
        the restrictions of universal
        families:
        \[
            \begin{tikzcd}
                \mathcal{M} \ar[d, no head] & \ar[l, hook']
                \mathcal{M}|_{\Spec\left(\OOhat_{\Quot_d(\cE), p_{C}}\right)}\ar[d, no head]
                \ar[r, "i", "\simeq"'] & \mathcal{M}'|_{\Spec\left(\OOhat_{
            \Quot_d(\cE'), p_{\mathbb{A}^1}}\right)}\ar[r, hook]\ar[d, no head] & \mathcal{M}' \ar[d, no head]\\
                \Quot_d(\cE) & \ar[l, hook'] \Spec(\OOhat_{
            \Quot_d(\cE), p_{C}}) \ar[r, "\iota", "\simeq"']
            &\Spec(\OOhat_{
        \Quot_d(\cE'), p_{\mathbb{A}^1}})\ar[r, hook] & \Quot_d(\cE')
            \end{tikzcd}
        \]
        Let $[q_{\bullet}']$ be a compatible quadric on $M'$ obtained from
        $[q_{\bullet}]$ using the isomorphism $\CQ(M') \simeq \CQ(M)$ given by
        $\iota$. The isomorphism $i$ yields an isomorphism
        \[
            \CCQ\left( \mathcal{M}|_{\Spec\left(\OOhat_{
        \Quot_d(\cE), p_C}\right)} \right)
            \simeq \CCQ\left( \mathcal{M}'|_{\Spec\left(\OOhat_{
            \Quot_d(\cE'), p_{\mathbb{A}^1}}\right)} \right),
        \]
        which shows that the complete local rings as $(M, [q_{\bullet}])$ and
        $(M', [q_{\bullet}]')$ are isomorphic.  Hence, we can and do reduce to the
        case $C = \mathbb{A}^1$, $\cE = \OO_{\mathbb{A}^1}^{\oplus r}$, in
        particular we drop the primes below.

        \def\CQU{\CQ(\mathcal{M}')}
        \def\Mtw{\widetilde{M}}
        \def\qtw{\widetilde{q}}
        We will use the infinitesimal lifting criterion. Let $\Spec(A)\into
        \Spec(B)$ be a small extension, where $B\onto A$ are Artinian local
        $\kk$-algebras. We are to prove that any deformation $(\Mtw,
        [\qtw_{\bullet}])$ of $(M, [q_{\bullet}])$ over $A$ can be
        extended to a deformation over $B$.

        The $A$-module $\Mtw$ is free. Take a free
        $B$-module $\Mtw'$ with surjection $\Mtw'\onto \Mtw$.
        The element $[\qtw_{\bullet}]$ is a map $\Spec(A)\to \CQ(\Mtw)\into
        \CQ(\Mtw')$.  The variety of completed quadrics $\CQ(\Mtw')$ is smooth
        over $B$, hence we take a lift $[q_{\bullet}]$ to a map
        $\Spec(B)\to \CQ(\Mtw')$ which corresponds to a broken quadric
        $[\qtw_{\bullet}']$.

        Let $\comp \to \CQ(\Mtw')$ be the compatible bundle.
        The broken quadric $[q_{\bullet}]$ is
        compatible with the action of $S = H^0(\OO_{\mathbb{A}^1}) = \kk[x]$, we have a structure morphism
        \[
            S\times \Spec(A) \to \comp|_{\Spec(A)} \to \CQ(\Mtw')|_{\Spec(A)}
        \]
        of $A$-schemes.
        Restricting it to the action of $x\in S$, we obtain a map
        \begin{equation}\label{eq:tmpComp}
            \Spec(A) \to \comp|_{\Spec(A)} \to \CQ(\Mtw')|_{\Spec(A)}
        \end{equation}
        The vector bundle
        $\comp$ is smooth over $\CQ(\Mtw')$, hence~\eqref{eq:tmpComp} lifts to
        a
        morphism $\Spec(B)\to
        \comp$, which yields a coherent $S_B$-module structure on $\Mtw'$, compatible with the broken quadric
        $[\qtw_{\bullet}']$. Since $\cE$ is free, the surjection $\cE|_{\Spec(A)}\onto \Mtw$ lifts
        to a surjection $\cE|_{\Spec(B)}\onto \Mtw'$ and so we obtain a map
        $\Spec(B)\to \CQuot_d(\cE)$, as required. This concludes the proof of
        smoothness.

        To prove connectedness of $\CQuot_d(\cE)$, consider first the open locus $V \subseteq
        \Quot_d(\cE)$ that parameterises modules $\cE\onto M$, where $M$ is
        supported at exactly $d$ points. The $\OO_C$-modules $M$ and
        $\OO_C/\!\Ann(M)$ are isomorphic. This implies that the fibre of
        $\CQuot_d(\cE)\to \Quot_d(\cE)$ over $[\cE\onto M]$ is isomorphic to the fibre of
        $\Iar_d(C)\to \OHilb_d(C)$ over $V(\Ann(M))\subseteq C$, hence it is
        isomorphism to the permutohedral variety, and so it is
        irreducible. Hence, the open
        locus $V$ is irreducible itself.

        To prove connectedness, it is enough
        to show that $V$ is dense. This amounts to producing, for every point,
        a germ $\Spec(\kkt)\to \Quot_d(\cE)$ such that the special point maps
        to a given point, while the generic points maps to $V$. As in the proof of
        smoothness, we reduce to the case
        $C = \mathbb{A}^1$, $\cE = \OO_{\mathbb{A}^1}^{\oplus r}$. At this
        point $V$ is no longer relevant for the proof.

        We will link any $\kk$-point $(M_0,
        [q_{\bullet}])$ to a special $\kk$-point which is $M_0 = \kk[x]/(x^d)$
        with an unbroken quadric coming from the functional $(x^{d-1})^*$
        as in Example~\ref{ex:quadricsOnGorenstein}.

        For the first step, consider the $\Gmult$-action on
        $\Quot_d(\OO_{\mathbb{A}^1}^{\oplus r})$, which is the composition of
        the usual action on $\mathbb{A}^1$ and an action with general weights
        assigned to $r$
        generators of $\OO_{\mathbb{A}^1}^{\oplus r}$. The only torus-fixed
        points are quotients of the form
        \begin{equation}\label{eq:tmpTorusFixed}
            \bigoplus_{i=1}^r \kk[x]e_i\onto \bigoplus_{i=1}^r
            \frac{\kk[x]}{(x^{\nu_i})} e_i.
        \end{equation}
        For the second step, observe that the
        infinitesimal lifting criterion above gives more than just smoothness:
        it shows that for any infinitesimal deformation of $q_{\bullet}$, the
        remaining structure deforms as well. This in particular implies that
        any $(M, [q_{\bullet}])$ lies in the closure of the locus $(M, [q_0])$
        where $q_{\bullet} = [q_0]$ is unbroken. This locus itself is
        irreducible, because it is open in the vector subspace
        $\Sym^2_A M\subsetneq \Sym^2_{\kk} M$, which is the symmetric
            square of \emph{an $A$-module, not a $\kk$-module}.

        For the third step, it remains to connect any module of the
        form~\eqref{eq:tmpTorusFixed} with a full rank quadric given by the
        functional $\sum_{i=1}^{r} (x^{\nu_i-1} e_i)^*$ to the module $M_0$
        above, which is a special case, with $0 = \nu_2 =  \ldots = \nu_r$.
        This is an elementary construction, which we leave to the reader. This
        concludes the proof of connectedness.
        The dimension estimate follows, since this is the dimension of the
        open locus $V$.
    \end{proof}

    \begin{theorem}\label{ref:smoothnessCurve:thm}
        Let $C$ be a smooth connected curve.
        The scheme $\Iar_d(C)$ is smooth irreducible of dimension
        $2d-1$. The morphism $\tau_C\colon \Iar_d(C)\to \OHilb_d(C)$
        is flat and projective with integral fibres of dimension $d-1$ which
        are reduced complete
        intersections in $\Iar_d(C)$. It satisfies $(\tau_C)_* \OO_{\Iar_d(C)}
        = \OO_{\OHilb_d(C)}$.
    \end{theorem}
    \begin{proof}
        Smoothness, irreducibility, and $\dim \Iar_d(C) = 2d-1$ all follow from
        Theorem~\ref{ref:smoothnessCQuotCurve:thm} for $\cE := \OO_C$.
        Projectivity of $\tau_C$ is clear, since $\tau_C$ is a restriction of
        $\tau_{\overline{C}}$ for a projective curve $\overline{C}$.
        We will show below that the fibres of $\tau_C$ are reduced, irreducible of dimension $d-1$.
        Once this is done, Miracle Flatness shows that $\tau_C$ is flat and that
        the fibres are complete intersections. Stein factorisation yields $(\tau_C)_* \OO_{\Iar_d(C)}
        = \OO_{\OHilb_d(C)}$.

        Take a subscheme $Z\subseteq C$ and embed it into $\mathbb{A}^1$. The
        fibres $\tau_C^{-1}(Z\subseteq C)$ and $\tau_{\mathbb{A}^1}^{-1}(Z
        \subseteq \mathbb{A}^1)$ are isomorphic, hence we reduce to the case
        $C = \mathbb{A}^1$.
        By Example~\ref{ex:completedQuadricsOnJordan} the fibre over
        $\Spec(\kk[x]/x^d)$ has dimension $d-1$.
        The morphism $\tau_{\mathbb{A}^1}$ is equivariant with respect
        to the $\Gmult$-action on $\mathbb{A}^1$.
        Every fiber of $\tau_{\mathbb{A}^1}$ is cut out by $d$ equations,
        hence has dimension $\geq 2d-1-d$ at every point. By semicontinuity
        using the torus action, every fiber has dimension $d-1$ everywhere and
        is a complete intersection.
        Finally, to prove reducedness of fibres, recall that they are
        Cohen-Macaulay,
        so~\cite[\href{https://stacks.math.columbia.edu/tag/031R}{Tag
        031R}]{stacks_project} it is enough to prove generic reducedness. But
        each of the fibres is over a Gorenstein subscheme, hence admits a
        point corresponding to an unbroken quadric, see
        Lemma~\ref{ref:unbrokenopen:lem}, which is smooth.
    \end{proof}

    \begin{example}
        The fibre of $\Iar_3(\mathbb{A}^1) \to \OHilb_3(\mathbb{A}^1)$ over
        $\Spec(\kk[x]/x^3)$ is singular. Indeed, by definition, the fibre is the locus in
        $\CQ(\kk[x]/x^3)$ of broken quadrics compatible with the action of
        $x$. We identify $R = \kk[x]/x^3$ with its dual $R^{\vee}$. One
        compatible quadric is the one corresponding to the full flag
        \[
            R = \frac{\kk[x]}{(x^3)} \supseteq (x)\supseteq (x^2)\supseteq (0).
        \]
        This quadric can be deformed according to the description from
        Example~\ref{ex:completedQuadricsOnJordan}. Namely, we can consider
        the quadrics with flag $R \supseteq (x^2)\supseteq (0)$ and a full
        rank quadric on $R/(x^2) \simeq \kk[x]/(x^2)$ given by
        \[
            \frac{R}{(x^2)}\times \frac{R}{(x^2)} \to \frac{R}{(x^2)}\to \kk,
        \]
        where at the end we apply an arbitrary functional that is nonzero on $x$. This yields one
        tangent vector. Performing an analogous procedure for the flag $R
        \supseteq (x)\supseteq (0)$ yields a mirrored situation and another
        tangent vector.

        The third tangent vector is not so easy to describe and it is best to
        employ Tyrrell's coordinates, as in~\S\ref{ssec:Tyrrell}. Namely,
        consider the ring $\kk[\varepsilon]/(\varepsilon^2)$ and the broken
        quadric which in Tyrrell's coordinates corresponds to $y_1 =
        \varepsilon$, $y_2 = -\varepsilon$, $x_{13} = \varepsilon$ and all
        other coordinates zero.
        In the exterior forms notation from~\S\ref{ssec:CQconstruction} this
        broken quadric is given by
        \[
             \begin{pmatrix}
                0 & 0 & \varepsilon\\
                0 & \varepsilon & 0\\
                \varepsilon & 0 & 0
            \end{pmatrix}\in \mathbb{P}\Sym^2 V,\qquad
            \begin{pmatrix}
                0 & 0 & -\varepsilon\\
                0 & -\varepsilon & 0\\
                -\varepsilon & 0 & 0
            \end{pmatrix}\in \mathbb{P}\Sym^2\Lambda^2 V
        \]
        This is compatible with $x$ and yields the third tangent vector. They can be checked to be linearly independent.
    \end{example}

    \begin{example}
        The morphism $\tau_{\Quot}\colon \CQuot_d(\OO_{\mathbb{A}^1}^{\oplus r})\to
        \Quot_d(\OO_{\mathbb{A}^1}^{\oplus r})$ is not flat for $r\geq 2$. Indeed, for
        example take the module $M = \kk[x]/(x) \oplus \kk[x]/(x^{d-1})$. Take
        any surjection that makes $M$ a point of the Quot scheme, and
        consider full rank compatible quadrics on $M$. The space of such quadrics is
        open in $\Sym^2_{\kk[x]}M\subsetneq \Sym^2 M$, which is a vector subspace of dimension
        $d+1$. We consider quadrics up to rescaling, but still this implies
        that the fibre of $\tau_{\Quot}$ is at least $d$-dimensional, while
        the fibre over a cyclic module $M$ is $(d-1)$-dimensional. The same argument shows that
        the map $\CQuot_d(\OO_{\mathbb{A}^1}^{\oplus r})\to \Sym^d \mathbb{A}^1$ is not flat as well:
        indeed, over a $d$-tuple of points, there is only one possible module and it is cyclic, hence
        the fibre over such a $d$-tuple will have dimension strictly smaller than the fibre over $d[0]$.
    \end{example}

    \subsubsection{Applications to characteristic numbers}\label{ssec:charNumbers}

        In this section we reprove one of the results
        of~\cite{Jagiella_Pielasa_Shatshila}, as explained
        in~\S\ref{ssec:junehuh}.

    \newcommand{\cL}{\mathcal{L}}%
        For any rank $d$ locally free sheaf $\cV$ on a scheme $B$,
        the embedding of $\CQ(\cV)$ into the product $\mathcal{P}$ given
        in~\eqref{eq:product} equips $\CQ(\cV)$ with $d-1$ line bundles
        $\cL_1, \ldots ,\cL_{d-1}$: the
        pullbacks of $\OO(1)$ from the first, second,  \ldots ,
        $(d-1)$-th factor.
        \begin{proposition}\label{ref:characteristicNumbers:prop}
            Let $\kk$ have characteristic different from two (we just repeat
            here the global assumption to stress that we do not assume
            characteristic zero).
            Let $b_1, \ldots ,b_{d-1}$ be any integers summing up to $d-1$.
            Let $V = \kk[x]/(f)$ be a finite $\kk$-algebra with
            $\dim_{\kk} V = d$.
            Let $X_{V}\subseteq \CQ(V)$ be the associated variety, as
            described in~\S\ref{ssec:junehuh}. Then the intersection number
            \[
                L_1^{b_1} \ldots L_{d-1}^{b_{d-1}}|_{X_V}
            \]
            is independent of $f$, that is, constant on the whole Hilbert
            scheme $\OHilb_d(\mathbb{A}^1)$.
        \end{proposition}
        \begin{proof}
            We will prove that the variety $X_V$ is the fibre of the
            Iarrobino-to-Hilbert map
            \[
                \tau_{\mathbb{A}^1}\colon \Iar_d(\mathbb{A}^1)\to
                \OHilb_d(\mathbb{A}^1).
            \]
            The variety $X_V$ is obtained as follows. The multiplication on
            $V$ is a map $\Sym^2 V\to V$. Its transpose is a map $V^{\vee} \to \Sym^2
            V^{\vee}$. We consider elements $\kk\to V^{\vee} \to \Sym^2
            V^{\vee}$ of its image which are full rank quadrics and we take the
            closure of this locus in $\CQ(V^{\vee})$.
            Observe that $\kk\to V^{\vee} \to \Sym^2 V^{\vee}$ is of full rank
            if and only if the dual map $\Sym^2 V\to V \to \kk$ is of full
            rank. By Example~\ref{ex:quadricsOnGorenstein} the locus above
            consists exactly of $q^{-1}$ for a full rank compatible quadric
            $q\in \Sym^2 V$.

            Recall from~\S\ref{ssec:dualityOnCQ} that there is an isomorphism
            $\CQ(V^{\vee})\to \CQ(V)$ which maps a full rank quadric $q$ to
            its inverse $q^{-1}$ and interchanges $\cL_{j}$ with
            $\cL_{d-j}$. Taking the closures of the corresponding
            loci, we obtain an isomorphism of $X_V$ and
            $\tau^{-1}_{\mathbb{A}^1}([f])$.

            Now, the desired characteristic number is the intersection number
            \[
                \cL_{d-1}^{b_1} \ldots
                \cL_{1}^{b_{d-1}}|_{\tau^{-1}_{\mathbb{A}^1}([f])}
            \]
            for a $\kk$-point $[f]\in \OHilb_d(\mathbb{A}^1)$. Since
            $\tau_{\mathbb{A}^1}$ is flat, this is independent of $[f]$,
            see~\cite[\S10.2]{Fulton__Intersection_theory}, since
            the intersection number is obtained from Euler
            characteristics.
        \end{proof}

        \subsubsection{Commuting matrices}\label{ssec:commuting}

    It is well-known that $\OHilb_d(\mathbb{A}^n)$ can be presented as a
    global quotient of the space of commuting matrices together with a
    generating vector and stability condition,
    see~\cite{nakajima_lectures_on_Hilbert_schemes}.
    As explained in \S\ref{ssec:introADHM}, the same holds for $\Iar_d(\mathbb{A}^n)$. In this section we prove the
    corresponding statements.

    \begin{proof}[Proof of Theorem~\ref{ref:ADHMdesciption:thm}]
        Fix a $d$-dimensional $\kk$-vector space $V$.
        Consider the universal rank $d$ vector bundle $\mathcal{U}$ on $\OHilb_d(\mathbb{A}^n)$.
        Let $H := \Isom(\mathcal{U}, V)\to \OHilb_d(\mathbb{A}^n)$ be the framing bundle, so that a morphism $T\to H$
        consists of a map $f\colon T\to \OHilb_d(\mathbb{A}^n)$ and a fixed isomorphism $f^*\mathcal{U} \to \OO_T\otimes_{\kk} V$.
        The morphism $f_H\colon H\to \OHilb_d(\mathbb{A}^n)$ is a $\GL(V)$-bundle and $f_H^*\mathcal{U}$ has a canonical trivialisation
        $f_H^*\mathcal{U}\simeq \OO_H\otimes_{\kk} V$.

        Consider the pullback $H' := H\times_{\OHilb_d(\mathbb{A}^n)} \Iar_d(\mathbb{A}^n)$. The Iarrobino scheme is a quotient of this pullback by $\GL(V)$.
        It remains to prove that $H'$ is isomorphic to the scheme given in the statement of Theorem~\ref{ref:ADHMdesciption:thm}.
        By Definition~\ref{ref:Iar:def} the pullback $H'$ is the compatible locus $\CCQ(f_H^*\mathcal{U})$. But $f_H^*\mathcal{U}\simeq \OO_{H}\otimes_{\kk} V$,
        so that $\CCQ(f_H^*\mathcal{U})\simeq \CCQ(V\times H')$ canonically. In
        particular, we obtain a map $H'\to \CQ(V\times H')\simeq \CQ(V) \times H'\to \CQ(V)$. By
        Definition~\ref{ref:compatibleLocus:def} of the compatible locus, we obtain
        $n$ commuting sections of the compatible bundle: the images of
        multiplications by variables in $H^0(\OO_{\mathbb{A}^n})$. The unit
        element $\OO_{\OHilb_d(\mathbb{A}^n)}\to \mathcal{U}$ yields the
        generating vector $v$. In this way we obtain a map from $H'$ to the scheme from the statement. The procedure
        can be reversed.
    \end{proof}
    In particular, for $n=1$, there are no commutativity conditions and so $\Iar_d(\mathbb{A}^1)$ is smooth
    of dimension $\binom{d+1}{2}-1 + \binom{d+1}{2} + d - d^2 = 2d-1$.

    The proof of Proposition~\ref{ref:ADHMGorenstein:prop} is similar and we leave it to the reader, only remarking that
    in the first description, the quadric is obtained from the usual scalar product.

\subsection{Torus limits on $\mathbb{A}^n$ and Iarrobino's symmetric
decomposition}\label{ssec:IarrobinoCurves}

In this section we explain the origin of the name Iarrobino scheme, based on
Anthony's Iarrobino work recalled in~\S\ref{ssec:symmDec}.

We consider $\Iar_d(\mathbb{A}^n)$ and fix the action of the torus $\Gmult$ on
$\mathbb{A}^n$ by scalar multiplication:
\[
    t\cdot (x_1, \ldots ,x_n) = (tx_1, \ldots ,tx_n).
\]
We will show that the limit at $t\to \infty$ of a Gorenstein algebra is a self-dual module described in~\S\ref{ssec:symmDec}.
This allows to link Iarrobino's symmetric decomposition to \BBname{} decompositions, which were successfully employed
to study the Hilbert scheme~\cite{Jelisiejew__Elementary, Jelisiejew__Pathologies}.

For a point $[Z]$ in the Hilbert scheme $\OHilb_d(\mathbb{A}^n)$ it is well-known that the limit 
$\lim_{t\to \infty} t\cdot [Z]$ exists only if $Z$ is supported exclusively at the origin. If this is the case,
then $Z = \Spec(A)$ for a local algebra $A$ and $\lim_{t\to \infty} t\cdot [Z]$ is isomorphic, as an abstract scheme,
to $\Spec(\gr(A))$.
The main point of this section is to derive an analogue of this observation for $\Iar_d(\mathbb{A}^n)$.

Consider a finite Gorenstein subscheme $Z = \Spec(A)\subseteq \mathbb{A}^n$ of degree $d$
supported only at the origin,
and fix an orientation $q\colon \omega_Z\to \OO_Z$ as in
Example~\ref{ex:quadricsOnGorenstein}. The pair $([Z], [q])$ is a $\kk$-point of
the unbroken locus $\Iar_d^{\Gor}(\mathbb{A}^n)$ from~\S\ref{ssec:unbroken}.
Since the limit $\lim_{t\to \infty} t\cdot [Z]$ in $\HilbdAn$ exists, and $\Iar_d(\mathbb{A}^n)\to \OHilb_d(\mathbb{A}^n)$ is
projective, also the limit
\[
    \lim_{t\to \infty} t\cdot ([Z], [q])
\]
exists in $\Iar_d(\mathbb{A}^n)$, say it is equal to $([Z_0], [q_{\bullet}])$, where $Z_0\simeq \Spec(\gr(A))$. It remains to show that $[q_{\bullet}]$
equips this algebra exactly with Iarrobino's symmetric decomposition.

\begin{example}\label{ex:associatedGradedTorusSymmetric}
    Let us verify the claim in the situation from
    Example~\ref{ex:associatedGradedSymmetric}.

    Recall that we consider $Z = V(xy, x^2 - y^3)\subseteq \mathbb{A}^2$. The
    corresponding algebra $A = \frac{\kk[x, y]}{(xy, x^2 - y^3)}$ is a
    complete intersection, hence it is Gorenstein. The algebra has a basis $1,
    x, y, y^2, y^3$. Let $\alpha\colon A\to \kk$ be the functional dual to
    $y^3$ in this basis and let $\varphi\in \Sym^2 A^{\vee}$ be the associated quadric:
    $\varphi(a, a') := \alpha(aa')$. It has full rank and as a map
    $\varphi\colon A\to A^{\vee}$ it is given by
    \[
        \begin{array}[h]{c c c c c c}
            a & 1 & y & y^2 & y^3 & x\\
            \midrule
            \varphi(a) & (y^3)^* & (y^2)^* & y^* & 1^* & x^*
        \end{array}
    \]
    We need to compute the dual quadric $q\in \Sym^2 A$. As a map $q\colon
    A^{\vee}\to A$, it is equal to $\varphi^{-1}$, so it is given by
    \[
        \begin{array}[h]{c c c c c c}
            a^* & 1^* & y^* & (y^2)^* & (y^3)^* & x^*\\
            \midrule
            q(a^*) & y^3 & y^2 & y & 1 & x
        \end{array}
    \]
    Consider now the $\Gmult$-orbit of $[Z]$. We will be interested in the
    limit at $t\to \infty$, so we use the coordinate $u := t^{-1}$. The action
    on $Z$ is given by
    \[
        u_0\cdot [Z] = \Spec\left( \frac{\kk[x, y]}{(xy, x^2 - u_0 y^3)} \right)
    \]
    for $u_0\in \Gmult(\kk)$.
    We see that for every nonzero $u_0$ the resulting quotient algebra
    \[
        A_{u_0} := \frac{\kk[x, y]}{(xy, x^2 - u_0y^3)}
    \]
    is again
    spanned by $1, y, y^2, y^3, x$ and equipped with a functional
    $\alpha_{u_0}$ defined by $\alpha_{u_0}(a) := u_0^3\cdot \alpha(u_0^{-1}a)$ for every
    $a\in A_{u_0}$. It follows that $\alpha_{u_0}$ is dual to $y^3$ in the
    above basis. From $\alpha_{u_0}$, we obtain $\varphi_{u}$ and
    $q_{u}$, given by
    \begin{align*}
        &\begin{array}[h]{c c c c c c}
            a & 1 & y & y^2 & y^3 & x\\
            \midrule
            \varphi_u(a) & (y^3)^* & (y^2)^* & y^* & 1^* & ux^*
        \end{array}\\
        &\begin{array}[h]{c c c c c c}
            a^* & 1^* & y^* & (y^2)^* & (y^3)^* & x^*\\
            \midrule
            q_u(a^*) & y^3 & y^2 & y & 1 & u^{-1} x
        \end{array}
    \end{align*}
    To compute $\lim_{t\to \infty} t\cdot ([Z], [q])$, it remains to compute
    $\lim_{u \to 0} q_u$. This is most conveniently done by
    applying~\S\ref{ssec:pointsInExterior}
    to $u\cdot q_u$, which yields the
    same class in $\mathbb{P}\Sym^2 A_u$ as $q_u$. It follows that $\lim_{u\to 0} q_u$ is the
    broken quadric corresponding to the flag
    \[
        \spann{1^*, y^*, (y^2)^*, (y^3)^*, x^*} \supseteq \spann{1^*, y^*,
        (y^2)^*, (y^3)^*} \supseteq (0)
    \]
    dual to the flag
    \[
        \frac{\kk[x, y]}{(x^2, xy, y^4)} \supseteq (x) \supseteq (0),
    \]
    and the (unique up to scalars) quadric $q_0 := (u q_u)|_{u=0}$ on $\spann{x^*}  \simeq \kk
    x$ together with the residual quadric $q_1$ on $\frac{\kk[x, y]}{(x^2, xy,
    y^4, x)}$ given by the functional $(y^3)^*$.
    This is in accordance with
    Example~\ref{ex:associatedGradedSymmetric}: the quadric $q_1$ yields
    self-duality of $\QHilb{0} \simeq \frac{\kk[x, y]}{(x, y^4)}$ and the
    quadric $q_0$ does it for $\kk x  \simeq \frac{\kk[x, y]}{(x, y)}$.
\end{example}

\begin{theorem}[Iarrobino decomposition is a torus limit on the Iarrobino
    scheme]\label{ref:name:thm}
    Let $Z = \Spec(A) \subseteq \mathbb{A}^n$ be supported only at the origin, where
    $A$ is local with maximal ideal $\mm$ and $s = \max\{i\colon \mm^i\neq
    0\}$.

    Let $q\colon
    \omega_Z \to \OO_Z$ be an orientation.
    The limit $\lim_{t\to \infty} t\cdot [Z]$ is given by $([Z_{\bullet},
        [q_{\bullet}])$, where $Z_{\bullet}$ is a flag of subschemes
    \[
        \mathbb{A}^n \supseteq Z_0\supseteq Z_1\supseteq  \ldots
    \]
    with $Z_0  \simeq  \Spec(\gr(A))$, and for $\delta\geq 1$, the scheme $Z_\delta$ is given by the ideal
    \[
        I_{\delta} = \bigoplus_{i\geq 0} \frac{\mm^i \cap \annmm{s+1-i-\delta} +
    \mm^{i+1}}{\mm^{i+1}} \simeq \bigoplus_{i\geq 0} \frac{\mm^i \cap
        \annmm{s+1-i-\delta}}{\mm^{i+1}\cap \annmm{s+1-i-\delta}}
    \]
    and
    \[
        q_i\colon \left(\frac{I_{Z_{\delta+1}}}{I_{Z_\delta}}\right)^{\vee}\to
            \frac{I_{Z_{\delta+1}}}{I_{Z_\delta}},
    \]
    proves that the subquotient
    \begin{equation}\label{eq:tmpQsub}
        \frac{I_{Z_{\delta+1}}}{I_{Z_\delta}}  \simeq \bigoplus_{i\geq 0} \frac{\mm^i \cap
            \annmm{s+1-i-\delta}}{\mm^{i+1}\cap \annmm{s+1-i-\delta} + \mm^{i}\cap
            \annmm{s-i-\delta}} = \QHilb{\delta}
    \end{equation}
    is self-dual.
\end{theorem}

\begin{proof}
    \def\cZ{\mathcal{Z}}%
    Let $S = \kk[x_1, \ldots ,x_n]$ and let $S^{\vee}  \simeq
    \kk[x_1^*, \ldots ,x_n^*]$ be the graded dual vector space. We have $\OO_Z
    \simeq S/I$ for an ideal $I$. Since $Z$ is supported only at the origin,
    the ideal $I$ contains $(x_1, \ldots ,x_n)^d$, so we have $\omega_Z = \OO_Z^{\vee}\into
    S^{\vee}$. Let $\varphi := q^{-1}$ and $\alpha := \varphi(1_{\OO_Z}, -)\colon \OO_Z\to \kk$ be a
    functional, then $(q^{-1})(a_1, a_2) = \alpha(a_1a_2)$, see
    Example~\ref{ex:quadricsOnGorenstein} and $\alpha\in \omega_Z\subseteq
    S^{\vee}$.

    From the definition of $s$, we have $\alpha =
    \alpha_s +  \ldots + \alpha_0$.
    Let us see how $\varphi$ interacts with $\QHilb{\bullet}$. Fix any $i,
    \delta$. Take a basis of
    $(\QHilb{\delta})_i$ from~\eqref{eq:tmpQsub} and lift it
    to a
    subset $\mathcal{B}_{i, \delta}\subseteq S_{i}$ of homogeneous elements of $S$ of degree $i$.
    For $b\in \mathcal{B}_{i, \delta}$ we have
    \[
        \varphi(b) = \alpha_s(b\cdot (-)) + \alpha_{s-1}(b\cdot (-)) +
        \ldots\in \OO_Z^{\vee}\subseteq S^{\vee}.
    \]
    This is an element annihilated by $\mm^{s+1-i-\delta}$. Being annihilated $\mm^{s+1-i-\delta}$ is the
    same as being annihilated by $(x_1, \ldots ,x_n)^{s+1-i-\delta}$, which is
    the same as being an element of $S^{\vee}_{\leq
    s-i-\delta}$. By homogeneity, for each $j$, the element $\alpha_j(b\cdot (-))$ lies in
    $S^{\vee}_{j-i}$, hence for $j \geq s-\delta+1$ we obtain
    \[
        \alpha_j(b\cdot (-)) = 0.
    \]
    Let us denote
    \[
        W_{i, \delta} := \left\{ \alpha_{s-\delta}(b)\ |\ b\in
        \mathcal{B}_{i, \delta} \right\} \subseteq S^{\vee}_{s-\delta-i}.
    \]
    By construction, the images of elements in $\mathcal{B}_{i, \delta}$ in $\gr(A)$ are
    linearly independent modulo the space $\mm^{i+1}\cap \annmm{s+1-i-\delta} + \mm^{i}\cap
    \annmm{s-i-\delta}$. This translates to saying that $W_{i, \delta}$ are
    linearly independent modulo
    \[
        \spann{W_{i+j, \delta-j}\ |\ j=1, \ldots ,\delta}.
    \]

    Let us now understand how the torus action on $S$ and $S^{\vee}$ fits
    into the picture. For the torus action, we have $\deg
    x_{\bullet}^* = -1$ and $\deg(x_i) = 1$.
    Take the parameter $u := t^{-1}$ and consider the morphism $\Spec(\kk[u])\to \HilbdAn$ that yields the limit at
    $t\to \infty$ and the corresponding subscheme $\cZ \subseteq \mathbb{A}^n
    \times \Spec(\kk[u])$.
    The ideal $I(\cZ)\subseteq
    S[u]$ is the homogenisation of $I$ using $u$. Correspondingly, for $u\neq
    0$, the functional $\alpha_u\in \OO_{\cZ|_{u\neq 0}}^{\vee}$ that yields the orientation is the
    homogenisation of $\alpha$ with respect to $u$:
    \[
        \alpha_u = \alpha_s + u \alpha_{s-1} + \ldots + u^s\alpha_0,
    \]
    It is homogeneous of degree $-s$ and yields the quadric
    \[
        \varphi_u(a_1, a_2) := \alpha_u(a_1a_2)\colon \OO_{\cZ|_{u\neq 0}}\to\OO_{\cZ|_{u\neq 0}}^{\vee}
    \]
    and $q_u =
    \varphi_u\colon \OO_{\cZ|_{u\neq 0}}^{\vee}\to \OO_{\cZ|_{u\neq 0}}$, both homogeneous of degree
    $-s$.
    For any $i$, $\delta$, and $b\in \mathcal{B}_{i, \delta}$, it follows from the discussion above that
    \[
        \varphi_u(b) = u^{\delta} \alpha_{s-\delta}(b) + u^{\delta+1}( \ldots
        )
    \]
    and that $\alpha_{s-\delta}(b)\in W_{i, \delta}$ are linearly independent
    modulo $\spann{W_{i+j, \delta-j}\ |\ j=1, \ldots ,\delta}$.
    Let
    $\mathcal{B}_{\delta} := \bigcup_{i=0}^s \mathcal{B}_{i, \delta}$ and
    $W_{\delta} := \bigcup_{i=0}^s W_{i, \delta}$. In the bases
    $\bigcup_{\delta} \mathcal{B}_{\delta}$ and $\bigcup_{\delta} W_{\delta}$, the matrix
    of $\varphi_u$ becomes
    \begin{equation}\label{eq:matrix}
        \begin{pmatrix}
            \Phi_{0} & 0 & 0 &  \ldots & 0\\
            * & u\Phi_{1} & 0 &  \ldots &0\\
            * & * & u^2\Phi_{2} &  \ldots &0\\
             \ldots
             * & * & * &  \ldots &u^{s}\Phi_{s}\\
        \end{pmatrix}
    \end{equation}
    where $u^{\delta}\Phi_{\delta}$ is a symmetric matrix corresponding to the map
    $\mathcal{B}_{\delta}\to W_{\delta}$; it is divisible exactly by
    $u^{\delta}$. (It is well-known, but irrelevant here, that $\Phi_s = 0$,
    $\Phi_{s-1} = 0$ whenever $s\geq 2$.)

    By~\cite[Theorem~7.20]{Kleiman_Thorup}, we can compute
    the limit of~\eqref{eq:matrix}  even though in our bases the matrix is
    non symmetric. This is because we compute the limit in each
    $\mathbb{P}(\Lambda^{\bullet}\gr(A))^{\otimes 2}$ and obtain a result in
    $\mathbb{P}\Sym^2(\Lambda^{\bullet}\gr(A)) =
    \left(\mathbb{P}(\Lambda^{\bullet}\gr(A))^{\otimes 2}\right)^{\mathrm{sym}}$.
    Since the matrix~\eqref{eq:matrix} is lower-triangular, the computation is
    straightforward: the resulting broken quadric $[q_{\bullet}]$ is a chain
    of quadrics $[\Phi_{0}],\,[\Phi_{1}],\, \ldots $ as claimed. Finally, by
    Duality~\S\ref{ssec:dualityOnCQ}, the limiting quadric for $q_u$ is
    obtained from the limiting quadric for $\varphi_u = (q_u)^{-1}$.
\end{proof}

    \begin{example}\label{ex:brokenGorensteinNotPrincipal}
            In this example we show that the subquotients for a point of $\Iar_d(\mathbb{A}^3)$ need not be
            principal modules, even in the closure of the Gorenstein locus.
            We consider $d=4$ and
            \[
                \Spec\left( \frac{\kk[x, y, z]}{(x, y, z)^2} \right) \subseteq
                \mathbb{A}^3.
            \]
            Let us take a family of Gorenstein subschemes converging to it. We
            choose a family
            \[
                \Spec\left( \frac{\kk[x, y, z][t]}{(x^2,\, y^2,\, tz -
                xy,\, z^2,\, xz,\, yz)} \right)\subseteq \mathbb{A}^3 \times \mathbb{A}^1
            \]
            which is the passage to the associated graded algebra
            of the complete intersection $(x^2,\, y^2,\, z-xy)$.
            Let us choose the functional $\alpha:=z^*$ as an orientation on
            this algebra.  It the usual bases $V = \spann{1, x, y, z}$ and
            $V^{\vee} = \spann{1^*, x^*, y^*,
            z^*}$, we have, for $t\neq 0$, the following
            \[
                    \varphi_t = \begin{pmatrix}
                        0 & 0 & 0 & 1\\
                        0 & 0 & t^{-1} & 0\\
                        0 & t^{-1} & 0 & 0\\
                        1 & 0 & 0 & 0
                    \end{pmatrix}
            \]
            so we obtain from $\varphi_t$ the broken quadric in
            $\CQ(V^{\vee})$ with the associated flag
            \begin{equation}\label{eq:tmpGor}
                V\supsetneq \spann{y, z} \supsetneq 0,
            \end{equation}
            and quadrics $[1^*\cdot x^*]$, $[y^*\cdot z^*]$ on subquotients.
            Correspondingly, the broken quadric $\lim_{t\to 0} \varphi_t^{-1}$
            on $\CQ(V)$ has the associated flag
            \[
                V^{\vee} \supsetneq \spann{1^*, z^*} \supsetneq (0)
            \]
            and quadrics $[y\cdot z]$, $[1\cdot x]$ on subquotients.
            We see that $\spann{y, z}$ is not a principal ideal,
            hence~\eqref{eq:tmpGor} does not yield a broken Gorenstein algebra
            structure.
        \end{example}


\begin{thebibliography}{{Woj}24}

\bibitem[B{\u{a}}l25]{Balibanu}
Ana B{\u{a}}libanu.
\newblock Wonderful varieties with a view towards {P}oisson geometry.
\newblock In {\em Advances in {P}oisson geometry}, Adv. Courses Math. CRM Barcelona, pages 85--112. Birkh{\"a}user/Springer, Cham, [2025] \copyright 2025.

\bibitem[BB73]{BialynickiBirula__decomposition}
A.~Bia{\l{}}ynicki-Birula.
\newblock Some theorems on actions of algebraic groups.
\newblock {\em Ann. of Math. (2)}, 98:480--497, 1973.

\bibitem[BE77]{BuchsbaumEisenbudCodimThree}
David~A. Buchsbaum and David Eisenbud.
\newblock Algebra structures for finite free resolutions, and some structure theorems for ideals of codimension {$3$}.
\newblock {\em Amer. J. Math.}, 99(3):447--485, 1977.

\bibitem[Ber12]{Bertin__punctual_Hilbert_schemes}
Jos\'{e} Bertin.
\newblock The punctual {H}ilbert scheme: an introduction.
\newblock In {\em Geometric methods in representation theory. {I}}, volume~24 of {\em S\'{e}min. Congr.}, pages 1--102. Soc. Math. France, Paris, 2012.

\bibitem[Bri07]{Brion__loghomogeneous}
Michel Brion.
\newblock Log homogeneous varieties.
\newblock In {\em Proceedings of the {XVI}th {L}atin {A}merican {A}lgebra {C}olloquium ({S}panish)}, Bibl. Rev. Mat. Iberoamericana, pages 1--39. Rev. Mat. Iberoamericana, Madrid, 2007.

\bibitem[CJN15]{cjn13}
Gianfranco Casnati, Joachim Jelisiejew, and Roberto Notari.
\newblock Irreducibility of the {G}orenstein loci of {H}ilbert schemes via ray families.
\newblock {\em Algebra Number Theory}, 9(7):1525--1570, 2015.

\bibitem[CM26]{Conner_Michalek}
Austin Conner and Mateusz Micha{\l}ek.
\newblock Characteristic numbers and chromatic polynomial of a tensor.
\newblock {\em Forum Math. Sigma}, 14:Paper No. e108, 2026.

\bibitem[CN09]{casnati_notari_irreducibility_Gorenstein_degree_9}
Gianfranco Casnati and Roberto Notari.
\newblock On the {G}orenstein locus of some punctual {H}ilbert schemes.
\newblock {\em J. Pure Appl. Algebra}, 213(11):2055--2074, 2009.

\bibitem[DCP83]{Concini_Procesi__Complete_symmetric_varieties}
C.~De~Concini and C.~Procesi.
\newblock Complete symmetric varieties.
\newblock In {\em Invariant theory ({M}ontecatini, 1982)}, volume 996 of {\em Lecture Notes in Math.}, pages 1--44. Springer, Berlin, 1983.

\bibitem[DCS99]{DeConcini_Springer}
C.~De~Concini and T.~A. Springer.
\newblock Compactification of symmetric varieties.
\newblock {\em Transform. Groups}, 4(2-3):273--300, 1999.
\newblock Dedicated to the memory of Claude Chevalley.

\bibitem[Eis95]{eisenbud}
David Eisenbud.
\newblock {\em Commutative algebra}, volume 150 of {\em Graduate Texts in Mathematics}.
\newblock Springer-Verlag, New York, 1995.
\newblock With a view toward algebraic geometry.

\bibitem[ER12]{elias_rossi_short_Gorenstein}
Joan Elias and Maria~E. Rossi.
\newblock Isomorphism classes of short {G}orenstein local rings via {M}acaulay's inverse system.
\newblock {\em Trans. Amer. Math. Soc.}, 364(9):4589--4604, 2012.

\bibitem[Fog68]{fogarty}
John Fogarty.
\newblock Algebraic families on an algebraic surface.
\newblock {\em Amer. J. Math}, 90:511--521, 1968.

\bibitem[Ful98]{Fulton__Intersection_theory}
William Fulton.
\newblock {\em Intersection theory}, volume~2 of {\em Ergebnisse der Mathematik und ihrer Grenzgebiete. 3. Folge. A Series of Modern Surveys in Mathematics [Results in Mathematics and Related Areas. 3rd Series. A Series of Modern Surveys in Mathematics]}.
\newblock Springer-Verlag, Berlin, second edition, 1998.

\bibitem[GLS07]{Gustavsen_Laksov_Skjelnes__Elementary_explicit_const}
Trond~S. Gustavsen, Dan Laksov, and Roy~Mikael Skjelnes.
\newblock An elementary, explicit, proof of the existence of {H}ilbert schemes of points.
\newblock {\em J. Pure Appl. Algebra}, 210(3):705--720, 2007.

\bibitem[Hai02]{Haiman__exposition}
Mark Haiman.
\newblock Notes on {M}acdonald polynomials and the geometry of {H}ilbert schemes.
\newblock In {\em Symmetric functions 2001: surveys of developments and perspectives}, volume~74 of {\em NATO Sci. Ser. II Math. Phys. Chem.}, pages 1--64. Kluwer Acad. Publ., Dordrecht, 2002.

\bibitem[Hui06]{Huibregtse_construction}
Mark~E. Huibregtse.
\newblock An elementary construction of the multigraded {H}ilbert scheme of points.
\newblock {\em Pacific J. Math.}, 223(2):269--315, 2006.

\bibitem[Iar94]{iarrobino_associated_graded}
Anthony Iarrobino.
\newblock Associated graded algebra of a {G}orenstein {A}rtin algebra.
\newblock {\em Mem. Amer. Math. Soc.}, 107(514):viii+115, 1994.

\bibitem[IMM21]{Iarrobino_Macias_Marques}
Anthony Iarrobino and Pedro Macias~Marques.
\newblock Symmetric decomposition of the associated graded algebra of an {A}rtinian {G}orenstein algebra.
\newblock {\em J. Pure Appl. Algebra}, 225(3):Paper No. 106496, 49, 2021.

\bibitem[Jel19]{Jelisiejew__Elementary}
Joachim Jelisiejew.
\newblock Elementary components of {H}ilbert schemes of points.
\newblock {\em Journal of the {L}ondon {M}athematical {S}ociety}, 100(1):249--272, 2019.

\bibitem[Jel20]{Jelisiejew__Pathologies}
Joachim Jelisiejew.
\newblock Pathologies on the {H}ilbert scheme of points.
\newblock {\em Invent. Math.}, 220(2):581--610, 2020.

\bibitem[JMR23]{Jelisiejew_Masuti_Rossi}
Joachim Jelisiejew, Shreedevi~K. Masuti, and M.~E. Rossi.
\newblock On the {H}ilbert function of {A}rtinian local complete intersections of codimension three.
\newblock {\em J. Pure Appl. Algebra}, 227(7):Paper No. 107326, 21, 2023.

\bibitem[JPS25]{Jagiella_Pielasa_Shatshila}
Jakub Jagie{\l}{\l}a, Pawe{\l} Pielasa, and Anatoli Shatsila.
\newblock Characteristic numbers of algebras.
\newblock {\em arXiv:19587}, 2025.

\bibitem[J{\v{S}}22]{jelisiejew_sivic}
Joachim Jelisiejew and Klemen {\v{S}}ivic.
\newblock Components and singularities of {Q}uot schemes and varieties of commuting matrices.
\newblock {\em J. Reine Angew. Math.}, 788:129--187, 2022.

\bibitem[Kle80]{Kleiman__history}
Steven~L. Kleiman.
\newblock Chasles's enumerative theory of conics: a historical introduction.
\newblock In {\em Studies in algebraic geometry}, volume~20 of {\em MAA Stud. Math.}, pages 117--138. Math. Assoc. America, Washington, DC, 1980.

\bibitem[KMR98]{kleppe_roig_codimensionthreeGorenstein}
Jan~O. Kleppe and Rosa~M. Mir{\'o}-Roig.
\newblock The dimension of the {H}ilbert scheme of {G}orenstein codimension {$3$} subschemes.
\newblock {\em J. Pure Appl. Algebra}, 127(1):73--82, 1998.

\bibitem[Kun11]{Kunte__Gorenstein_modules_of_finite_length}
Michael Kunte.
\newblock Gorenstein modules of finite length.
\newblock {\em Math. Nachr.}, 284(7):899--919, 2011.

\bibitem[Lak87]{Laksov__completed_quadrics}
Dan Laksov.
\newblock Completed quadrics and linear maps.
\newblock In {\em Algebraic geometry, {B}owdoin, 1985 ({B}runswick, {M}aine, 1985)}, volume 46, Part 2 of {\em Proc. Sympos. Pure Math.}, pages 371--387. Amer. Math. Soc., Providence, RI, 1987.

\bibitem[Mic23]{Michalek__survey}
Mateusz Micha{\l}ek.
\newblock Enumerative geometry meets statistics, combinatorics, and topology.
\newblock {\em Notices Amer. Math. Soc.}, 70(4):588--597, 2023.

\bibitem[MS05]{Miller_Sturmfels}
Ezra Miller and Bernd Sturmfels.
\newblock {\em Combinatorial commutative algebra}, volume 227 of {\em Graduate Texts in Mathematics}.
\newblock Springer-Verlag, New York, 2005.

\bibitem[Nak99]{nakajima_lectures_on_Hilbert_schemes}
Hiraku Nakajima.
\newblock {\em Lectures on {H}ilbert schemes of points on surfaces}, volume~18 of {\em University Lecture Series}.
\newblock American Mathematical Society, Providence, RI, 1999.

\bibitem[Ols16]{Olsson}
Martin Olsson.
\newblock {\em Algebraic spaces and stacks}, volume~62 of {\em American Mathematical Society Colloquium Publications}.
\newblock American Mathematical Society, Providence, RI, 2016.

\bibitem[Pez18]{Pezzini__Lectures}
Guido Pezzini.
\newblock Lectures on wonderful varieties.
\newblock {\em Acta Math. Sin. (Engl. Ser.)}, 34(3):417--438, 2018.

\bibitem[Sem48]{Semple}
J.~G. Semple.
\newblock On complete quadrics.
\newblock {\em J. London Math. Soc.}, 23:258--267, 1948.

\bibitem[sta26]{stacks_project}
{S}tacks {P}roject.
\newblock \url{http://math.columbia.edu/algebraic_geometry/stacks-git}, 2026.

\bibitem[Tha99]{Thaddeus__complete_collineations}
Michael Thaddeus.
\newblock Complete collineations revisited.
\newblock {\em Math. Ann.}, 315(3):469--495, 1999.

\bibitem[TK88]{Kleiman_Thorup}
Anders Thorup and Steven Kleiman.
\newblock Complete bilinear forms.
\newblock In {\em Algebraic geometry ({S}undance, {UT}, 1986)}, volume 1311 of {\em Lecture Notes in Math.}, pages 253--320. Springer, Berlin, 1988.

\bibitem[Tyr56]{Tyrrell__complete_quadrics}
J.~A. Tyrrell.
\newblock Complete quadrics and collineations in {$S_n$}.
\newblock {\em Mathematika}, 3:69--79, 1956.

\bibitem[{Woj}24]{Wojtala}
Maciej {Wojtala}.
\newblock Iarrobino's decomposition for self-dual modules.
\newblock arXiv:2405.13829, 2024.

\end{thebibliography}
\end{document}